\documentclass[a4paper,12pt]{article}
\makeatletter
\let\ams@underbrace=\underbrace
\def\underbrace#1_#2{%
  \setbox0=\hbox{$\displaystyle#1$}%
  \ams@underbrace{#1}_{\parbox[t]{\the\wd0}{#2}}%
}
\makeatother
\usepackage{float}
\usepackage{pst-all}
\usepackage{pstricks}
\usepackage{mathtools}
\usepackage{appendix}
\usepackage{enumerate}
\usepackage{amssymb}
\usepackage{stmaryrd}
\usepackage[standard]{ntheorem}
\everymath{\displaystyle}

\newtheorem{Prop}{Proposition}
\newcommand{\cg}{\llbracket }
\newcommand{\cd}{\rrbracket }
\newcommand{\CG}{\biggr\llbracket}
\newcommand{\CD}{\biggr\rrbracket}

\newcommand{\F}{\mathbb{F}}
\newcommand{\Q}{\mathbb{Q}}

\newcommand{\PP}{\mathcal{P}}

\newcommand{\N}{\mathbb{N} }
\newcommand{\Z}{\mathbb{Z} }

\title{The local lifting problem for $(\Z/2\Z)^3$}
\author{Guillaume Pagot}
\date{January 2024}

\begin{document}
\maketitle

%\section{Introduction}
%\begin{center}
%    \Large{\bf{Relèvements d'actions de groupe $(\Z/2\Z)^3$}}
%\end{center}

\begin{abstract}
Let $k$ be an algebraically closed field of characteristic $2$. In this paper we describe the $(\Z/2\Z)^3$-actions on $k[[z]]$ for which there is a discrete valuation ring $R$, a finite extension of the ring of Witt vectors $W(k)$, such that they can be lifted as a group of $R$-automorphisms of $R[[Z]]$. In fact the necessary and sufficient condition for such an action to lift involves only the conductor type of the corresponding extension.

\end{abstract}

\section{Introduction}

%\section{Le résultat}
%\subsection{Notations et rappel du problème}

%  On considère ici un corps $k$ algébriquement clos de caractéristique $2$ et $R$ une extension "suffisamment grande" de l'anneau des vecteurs de Witt $W(k)$. L'anneau $R$ est un anneau de valuation discrète dont on notera $\pi$ une uniformisante. Notons $G$ le groupe $(\Z/2\Z)^3$. 
  
%  On se donne un revêtement $k[[z]]/k[[z]]^G$ de type $(m_1 + 1, m_2+1, m_3 + 1)$ relativement à des groupes $G_1,G_2,G_3$ (qui sont des sous-groupes d'indices $2$ de $G$) au sens de la définition 3.6 de Yang. Ce triplet $(m_1,m_2,m_3)$ est en particulier minimal.
  
  %On se donne une extension galoisienne $k[[z]]$ de $k[[t]]$ de groupe $G$. Celle-ci peut être donnée par %trois équations d'Artin-Schreier de conducteurs $m_1+1\leq m_2+1\leq m_3+1$, le triplet $(m_1,m_2,m_3)$ %étant choisi minimal.

 % On suppose seulement que $m_1+1$ est divisible par 4 ($m_2+1$ et $m_3+1$ devant bien sûr être pairs). On montre alors que l'on peur relever $G$ en un groupe d'automorphismes de $R[[Z]]$. On utilise le théorème 3.13 de ~\cite{Yang} qui nous dit qu'il "suffit" de relever chacun des trois revêtements intermédiaires d'ordre $2$ de façon à ce que les lieux de branchement de ces trois revêtements aient la "bonne combinatoire".  

%\subsection{The Theorem}
Let $k$ be an algebraically closed field of charasteristic $p>0$ and let $R$ be a complete discrete valuation ring, finite extension of the ring of Witt vectors $W(k)$. Let $\pi$ denote a uniformizing parameter of $R$ and $K=Fr(R)$. 

Let $n$ be a positive integer. Let $G$ be the group $(\Z/p\Z)^n$ and assume $G$ be a group of $k$-automorphisms of the ring $k[[z]]$ (in other words, we have an injection $G\hookrightarrow Aut_kk[[z]]$). We then write $k[[t]]=k[[z]]^G$, so that $k[[z]]$ is a Galois extension of $k[[t]]$ with Galois group $G$. Artin-Schreier theory tells us that this extension is the compositum of $n$ $p$-cyclic extensions of $k[[t]]$, given by equations of the form $y_i^2-y_i=f_i\left(\frac{1}{t}\right)$ for $i\in\cg1,n\cd$, where the $f_i$ are polynomials of degree $m_i$ prime to $p$ and the $f_i$ are $\F_p$ linearly independant for $1\leq i\leq n$. Reordering the equations, we can assume that $m_1\leq m_2\leq\cdots\leq m_n$ and that the $n$-tuple $(m_1,\cdots,m_n)$ is chosen minimally (for a more precise definition, we refer to definition 1 borrowed from Yang). The elements $m_1+1,\cdots,m_n+1$ are called the conductors of these intermediate extensions.

We are interested here in the problem of lifting $G$ as a subgroup of $Aut_R(R[[Z]])$.

Such a lifting is not always possible and obstructions are already known. In ~\cite{GrMa}, which deals with the case $n=2$, it is proved that if such lifting exists then necessarily $p$ divides $m_1+1$. For any $n$, this necessary condition generalizes and we must have the following divisibilities: for any $j\in\cg 1,n-1\cd$, $p^{n-j}|m_j+1$ (see ~\cite{Ber}). 

These obstructions are not the only ones : even if the divisibility conditions are met, there are obstructions of metric and differential nature as soon as $p>2$ and we hope that divisibility conditions are then the only possible obstructions when $p=2$. For further details on these obstructions, we refer to ~\cite{Pag} (Chapter 3).

This paper is a contribution to the study of a possible lifting in the case where $p=2$. This has already been proved in the case $n=2$ (see ~\cite{Pag} and  ~\cite{Mat2} ) and we give a slightly modified proof. The new result concerns the case $n=3$, and we prove in this case that under the only condition $4|m_1+1$, lifting is possible. More precisely , we prove the following theorem :

\begin{theorem}

Let $k$ an algebraically closed field of characteristic $2$. Let $n\in\N^*$ (the set of positive integers). Let $G$ be a group of $k$-automorphisms of the ring $k[[z]]$ and $G\simeq (\Z/2\Z)^n$. We then have the following two results:
    
    \begin{enumerate}[a.]
        \item If $n=2$, there exists a discrete valuation ring $R$,  which is a finite extension of the ring of Witt vectors $W(k)$, such that $G$ lifts to a group of $R$-automorphisms of $R[[Z]]$. 
         \item  If $n=3$, we always denote $(m_1+1,m_2+1,m_3+1)$ the "minimal triplet of conductors" for the $G$-action as group of $k$-automorphisms of $k[[z]]$ (if we use the terminology of Definition 1 below, this means that $k[[z]]/k[[z]]^{G}$ is a cover of type $(m_1+1,m_2+1,m_3+1)$). Then there exists a discrete valuation ring $R$, which is a finite extension of the ring of Witt vectors $W(k)$, such that $G$ lifts to a group of $R$-automorphisms of $R[[Z]]$ if and only if $4$ divides $m_1+1$. 
    \end{enumerate}
\end{theorem}

The present manuscript builds on notes
written in 2003-2004 after my thesis, which among other things dealt with the case of
the Klein four-group.  The present manuscript also draws on recent work
of Jianing Yang (see Theorem 2 below). On the sound advice of David Harbater and Michel Matignon, I have decided to write up these results and make them available to the mathematical community. I would like to take this opportunity to warmly thank Michel for his advice and proof-reading, and David for his encouragement and interest in this work.  

\section{The (local) lifting problem for curves}

Let $k$ be an algebraically closed field of charasteristic $p>0$ and let $R$ be a complete discrete valuation ring, finite extension of the ring of Witt vectors $W(k)$. Let $\pi$ denote a uniformizing parameter of $R$ and $K=Fr(R)$. 
 
Let $C$ be a smooth, projective curve over $k$ and $G$ be a finite subgroup of $Aut_k(C)$, so that $C \rightarrow D = C/G$ is a
branched $G$-cover of smooth, projective curves over $k$. We are interested in the following question : is there a ring $R$ (as described above) and a branched $G$-cover $ \mathcal{C} \rightarrow \mathcal{D}$ of smooth relative $R$-curves which lifts the cover $C \rightarrow D$ ? This question is known as "the lifting problem for covers of curves".

The answer to this question is already known in a few specific cases.

\begin{enumerate}[$\bullet$]
  \item If $|G|$ is not divisible by $p$, the answer is yes, by Grothendieck, SGA I.
  \item If $G$ is cyclic the answer is yes. This result, formerly known as the Oort conjecture, was jointly demonstrated by Obus-Wewers and
Pop (~\cite{OW14},  ~\cite{Pop14}).
  \item If $|G|>84(g(C)-1)$ with $g(C)\geq 2$, the answer is no. Indeed, in characteristic 0 the order of the automorphism group of a curve of genus $g\geq 2$ cannot exceed $84(g-1)$, whereas in characteristic $p$, there are curves with an automorphism group greater than $84(g-1)$ (~\cite{Ro}).
\end{enumerate}

In fact, a local-global principle (see for example ~\cite{Ga}, ~\cite{BM}, ~\cite{Hen} or~\cite{Ob19}) leads us  to consider the following question :

Suppose $G$ is a finite group, and assume $G$ be a group of $k$-automorphisms of the ring $k[[z]]$ (in other words, we have an inclusion $G\hookrightarrow Aut_kk[[z]]$). We can then write $k[[t]]=k[[z]]^G$, so that $k[[z]]$ is a Galois extension of $k[[t]]$ with Galois group $G$. Does there exist a discrete valuation ring $R$ of charasteristic 0 with residue field $k$ and a $G$-Galois extension $R[[Z]]/R[[T]]$ such that the $G$-action on $R[[Z]]$ gives the $G$-action on $k[[z]]$ in reduction ?

 We'll refer to this last question as" the local lifting problem". A group $G$ for which every local $G$-extension lifts to characteristic zero
is called a local Oort group for $k$. If there exists a local $G$-extension lifting to characteristic zero, $G$ is called a weak local Oort group.

Significant progress has been made in recent years concerning local Oort groups. We already know that the only possible local Oort groups are cyclic, dihedral of the form $D_{p^n}$, and $A_4$ (if $char(p)=2)$ (this a consequence of Theorem 1.2 in ~\cite{CGH}). More precisely, we know that the cyclic groups and $A_4$ for $p=2$ (~\cite{Ob16}) are local Oort groups. For dihedral groups, we know that $D_p$ is a local Oort-group (see ~\cite{BW} for $p>2$), as are $D_4$ (~\cite{Wea}), $D_9$ (~\cite{Ob17}), $D_{25}$ and $D_{27}$ (~\cite{DDKOT}). Note that Kontogeorgis and Terezakis have announced that some $D_{125}$ action does not lift (~\cite{KT}).

If a group $G$ is a weak Oort group and not an Oort group, we can look more closely and ask which $G$-covers lift. This is the question we propose to study for the groups $G=(\Z/p\Z)^2$ and  $G=(\Z/2\Z)^3$. In fact little is known for the actions of $G=(\Z/p\Z)^n$ with $n\geq 2$ and $p>2$ which can be lifted : as mentioned earlier, when $p>2$, differential constraints seriously complicate the situation.

\section{Tools for the proof}

We return to the notations used in I. Let's give definition 3.6 of ~\cite{Yang} again.

\begin{definition}

Let $G=(\Z/p\Z)^n$, and assume that $G$ is a group of $k$-automorphisms of
$k[[z]]$ as a $k$-algebra. Suppose that $(m_1, \cdots, m_n)$ is the lexicographically smallest n-tuple of
integers such that there exists subgroups $G_1,\cdots,G_n\subset G$ of index $p$ that satisfy the following
conditions :
\begin{enumerate}[1)]
\item The $p$-cyclic extensions $k[[z]]^{G_i}/k[[z]]^{G}$ have conductors $m_i + 1$,
\item  $k[[z]]^{G_1},\cdots ,k[[z]]^{G_i}$ are linearly disjoint over $k[[z]]^G$.
\end{enumerate}
Then we say that $k[[z]]/k[[z]]^{G}$ is a cover of type $(m_1+1,\cdots,m_n+1)$, with respect to  $(G_1,\cdots,G_n)$. Recall that $p$ does not divide any of the integers $m_i$ and note that $m_1 \leq m_2 \leq \cdots \leq m_n$.
\end{definition}

We already know of a few examples where lifting is possible when $G=(\Z/p\Z)^n$ (namely $G$ is a weak Oort local group). In order to prove this in ~\cite{Mat}, one gives examples of such liftings (with equidistant geometry) in the case of extensions of type $(p^{n-1}(p-1),\cdots,p^{n-1}(p-1))$. In ~\cite{Pag}, we show examples with non-equidistant geometry in the case of extensions of type $(qp^n,\cdots,qp^n)$ (for all $q\in\N^*$). More recently, Yang has given examples for $p=2$ and $n=3$, in the case of extensions of type $(4,4,2r)$ and $r$ is an integer greater than 1, which provides the first example where $n\geq 3$ and the $m_i$ are not all equal (see  ~\cite{Yang}).

A few reminders about order $p$-automorphisms of the $p$-adic open disk are necessary. Let $\sigma$ be an order $p$-automorphism of the $p$-adic open disk $Spec R[[Z]]$, which does not induce residually the identity. Let $F$ be the set of fixed points of $\sigma$ which can be assumed to be rational on $R$ (we may enlarge $R$ if necessary). Note that $|F|>1$. Consider the minimal semi-stable model of the marked open disk $(Spec R[[Z]],F)$. The special fiber is a tree of projective lines linked to the original generic point $(\pi)$ (which we will represent in the following figures in the form of a zigzag) by double crossing points. The fixed points specialize in the terminal components of the tree. In particular, each terminal component has at least two fixed points (this is an important consequence (prop 1.2. of ~\cite{GrMa2}) of the identity (*) below.

Each of the double points on the special fiber has a corresponding open annulus on the generic fiber (i.e. the points of the open annulus specialize into the double point). This annulus has a thickness, and we will indicate this thickness on the various figures to come with the letter "epsilon".

Let us briefly recall the result of Proposition 1.1 of \cite{GrMa2}.
Write $F=\{Z_0,Z_1,\cdots,Z_m\}$ and consider one of the fixed points $Z_i$. 
 Given $\rho\in R^{alg}$
with $v(\rho)\geq 0$, after enlarging $R$ so that $\rho\in R$, we let $v_{\rho}$ be the Gauss valuation
on $Fr(R[[Z]])$ relative to $\frac{Z-Z_i}{\rho}$
and $d(v(\rho))$ be the degree of the different of
the $v_\rho$-valued extension $Fr(R[[Z]])/ Fr(R[[Z]]^{<\sigma>})$. The graph of $d(v(\rho))$ for $v(\rho) \in \Q_+$ is piecewise linear with breaks in the set $\{v(Z_j-Z_i), j \neq i\} := \{v(\rho_1), v(\rho_2), \cdots, v(\rho_{\ell_i})\}$, where $v(\rho_1) < v(\rho_2) < ... < v(\rho_{\ell_i} )$.
For $k\in\cg1,\ell_i\cd$, let $\mu_k$ be the cardinality of the set of $Z_j \in  F$ such that $v(\rho_k) \leq v(Z_j-Z_i) <\infty$ and let $\mu_{\ell_i+1}=0$. We have then 
$$ \sum_{j=1}^{\ell_i}(\mu_j-\mu_{j+1})v(\rho_j)=\frac{1}{p-1}(v(p)-d(0)). \;\;\;(*)$$

Moreover, since $\sigma$ does not induce residually the identity, it follows that $d(0)=0$.

These geometric properties are formalized in the Hurwitz tree concept. Each projective line corresponds to a vertex, and each crossing point corresponds to an edge in the Hurwitz tree. Moreover, for each fixed point $Z_i$, we append a vertex and connect it via an edge to the vertex representing the projective line containing the specialization of $Z_i$. 
This Hurwitz tree is a metric tree (each edge corresponding to a double point is associated with the thickness of the corresponding annulus and the edges corresponding to fixed points are associated with 0). The preceding equality is then a metric condition that must be satisfied by any Hurwitz tree. 

In the following, we will focus exclusively on the branch locus so we will consider the Hurwitz tree associated to the minimal semi-stable model of the marked open disk $(Spec R[[T]],B)$ where  $R[[T]]=R[[Z]]^{<\sigma>}$ is the quotient of 
$(Spec R[[Z]],F)$. This new tree is the same Hurwitz tree except that the thicknesses are multiplied by $p$. When considering this quotient tree, equality (*) becomes 
$$ \sum_{j=1}^{\ell_i}(\mu_j-\mu_{j+1})v(\rho_j)=\frac{1}{p-1}pv(p). \;\;\;(**)$$

This equality will be crucial in part V.2. where we'll have to choose "good thicknesses". For further details and explanations, we refer to ~\cite{GrMa2} and ~\cite{Hen}.

When it's possible to lift this action, we can consider each of the $n$ covers lifting the extensions $k[[z]]^{G_i}/k[[z]]^{G}$ and consider their branch points, which will correspond to the fixed points of the order $p$ automorphism of the underlying $p$-adic open disk cover which lifts $k[[t]]^{G_i}$. We speak of equidistant geometry when the mutual distance between the fixed points in $Spec R[[Z]]^G$ is constant. Note that this is equivalent to the fact that these branch points specialize on the same terminal component.

A few clarifications concerning Theorem 1.
Assertion a. has already been proved for some twenty years (but not published in a journal and only available online). We propose to overcome this missing by including the proof (which differs slightly in calculation but not in spirit from that given in  ~\cite{Pag} and ~\cite{Mat2}) in the present article.

%A recent result by Yang (see Proposition 5.9 in ~\cite{Yang}) deals with case b. in the %special case where $m_1=m_2=4$. The present article is therefore a generalization of this %result.

In order to prove this theorem, we will use a recent result by Yang (see Theorem 3.13 in ~\cite{Yang}), generalizing an older result due to Green and Matignon, which gives a necessary and sufficient condition for an action of $G=(\Z/p\Z)^n$ on $k[[z]]$ to lift into an action on $R[[Z]]$. We give it here in its entirety.

\begin{theorem} \textbf{(Yang)}

Let $G=(\Z/p\Z)^n$. Suppose $k[[z]]/k[[t]]$ is a $G$-extension of conductor type $(m_1+1,\cdots,m_n+1)$ , with respect to $G_1,\cdots,G_n$. Then there is
a lift of $G$ to a group of automorphisms of $R[[Z]]$ if and only if the following two conditions hold:

\begin{enumerate}[a.]
    \item  $m_i\equiv -1 \mod p^{n-i}$ for $1 \leq i \leq n-1$,
    \item $k[[z]]^{G_1}/k[[z]]^G,\cdots k[[z]]^{G_n}/k[[z]]^G$ can be lifted with branch loci $B_1,\cdots,B_n$ such that for any subset
of $r$ branch loci $\{B_{i_1},\cdots,B_{i_r}\}$, $ |\bigcap_{1\leq j\leq r}B_{i_j}|=\frac{(\min_j(m_{i_j})+1)(p-1)^{r-1}}{p^{r-1}}$.
\end{enumerate}
   
\end{theorem}

The second condition is clearly of a combinatorial nature on the branch points and thus on Hurwitz trees of the $n$ Artin-Schreier coverings that induce the type.

So, for the proof of Theorem 1, we'll consider the two or three Artin-Schreier equations (for the cases $n=2$, $n=3$ respectively) giving the extension $k[[z]]/\!\raisebox{-.65ex}{\ensuremath{k[[z]]^G}}$. The idea is to find explicit equations of intermediate covers of the $p$-adic open disk, having simultaneously good reductions relative to the same Gauss valuation (these reductions of course giving the Artin-Schreier equations fixed beforehand). The tricky part is to ensure that the branch points of the covers thus constructed have the "right combinatorics", i.e. the one mentioned in the previous theorem. The proof is necessarily technical, as it involves giving explicit equations of the liftings.

Finally, to verify the good reduction of the intermediate $p$-order covers we'll be constructing, it will suffice (by virtue of result 3.4 in ~\cite{GrMa}) to check that the degrees of the generic and special differents are identical, which here amounts to establishing (since $p=2$) that the number of branch points is equal to the conductor of the corresponding Artin-Schreier extensions.

\section{The case $G=(\Z/2\Z)^2$}

The problem is already solved in the case $G=(\Z/2\Z)^2$ (see the unpublished manuscripts ~\cite{Pag} and ~\cite{Mat2}). We give a slightly different proof, which makes easier the reading of the proof in the case $G=(\Z/2\Z)^3$.

\subsection{Geometry of the branch locus}

   % Rappelons quelle était la géométrie du lieu de branchement.

   We start from an extension $k[[z]]/\!\raisebox{-.65ex}{\ensuremath{k[[z]]^G}}$, of type $(m_1+1,m_2+1)$, given by two Artin-Schreier equations :

$$\left\{ 
\begin{aligned}
    y_1^2-y_1&=f_1\left(\frac{1}{t}\right)\\
    y_2^2-y_2&=f_2\left(\frac{1}{t}\right)\\
\end{aligned}
\right.$$
where $f_1$ and $f_2$ are polynomials of degrees $m_1$ and $m_2$ respectively, with the condition $m_1\leq m_2$. 

The task is to find two covers of the $p$-adic open disk (denoted $Rev_1$ and $Rev_2$) which lift the two covers given above and such that $Rev_1$ and $Rev_2$ have exactly $\frac{m_1+1}{2}$ branch points in common (this is given by Theorem 2). Since we want to deal with the general case (and therefore in particular the case where $m_1<m_2$), we can't be satisfied with an equidistant geometry (i.e. a tree with a single projective line where the branch points would specialize). We therefore propose the geometry shown in Figure 1.

$ $
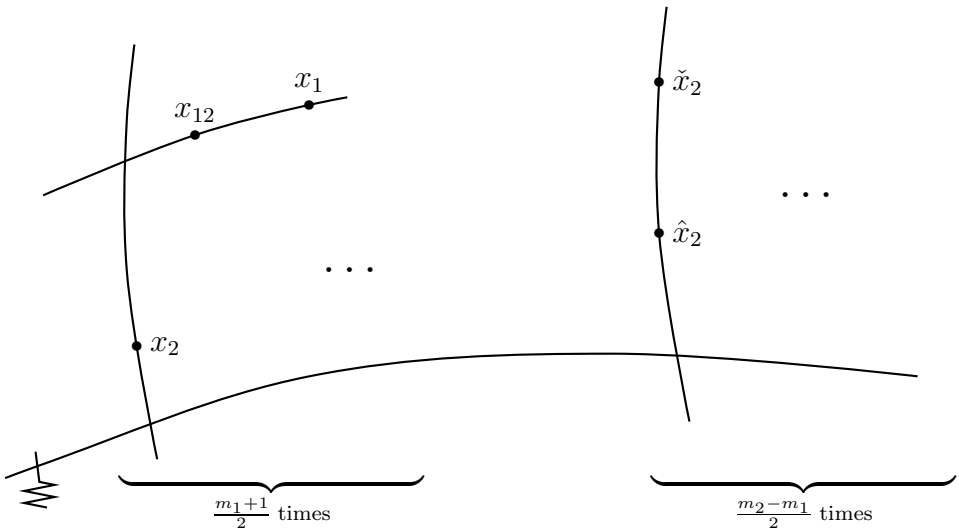
\begin{figure}[h]
\begin{pspicture}(0,-0.5)(16,6)	
%\pscurve[showpoints=true]{-}(0,1.3)
%(0.7,1.8)(3.3,0.5)(4,1.6)(0.4,0.4)
\psline[showpoints=false]{-}(0.4,0.6)(0.45,0.2)(0.65,0.16)(0.25,0.08)(0.65,0.0)(0.25,-0.08)(0.55,-0.14)%generique
\pscurve[showpoints=false]{-}(0,0.25)(0.4,0.4)(4,1.6)(8,1.9)(12,1.6)%horizontale0
\pscurve[showpoints=false]{-}(2,0.5)(1.6,3)(1.6,5)(1.7,6)%verticale1
\pscurve[showpoints=false]{-}(0.5,4)(2.5,4.8)(4,5.2)(4.5,5.3)%horizontale1
\pscurve[showpoints=false]{-}(9,1)(8.6,3.5)(8.6,5.5)(8.7,6.5)%verticale2

\psdot(1.73,2)\uput[0](1.73,2){$x_2$}
\psdot(2.5,4.8)\uput[90](2.5,4.8){$x_{12}$}
\psdot(4,5.2)\uput[90](4,5.2){$x_{1}$}
\psdot(8.6,3.5)\uput[0](8.6,3.5){$\hat{x}_{2}$}
\psdot(8.6,5.5)\uput[0](8.6,5.5){$\check{x}_{2}$}

\uput[270](3.5,0.5){$\underbrace{\hspace{4cm}}_{\frac{m_1+1}{2}\;\mathrm{times}}$}
\uput[270](10.5,0.5){$\underbrace{\hspace{4cm}}_{\frac{m_2-m_1}{2}\;\mathrm{times}}$}

\uput[0](4,3){\Large $\cdots$}\uput[0](10,4){\Large $\cdots$}
%\psaxes{->}(0,0)(-0.5,-0.5)(14.5,7.5)[$x$,0][$y$, 90]	%creates axes
%\psdot(2,1)	%plots the point (2,1)
%\uput[0](2,1){$A$}	%labels the point (2,1) as A

\end{pspicture}
\caption{Branch locus when $G=(\Z/2\Z)^2$}\label{fig1}

\end{figure}

A few comments on this figure are necessary, as the notations differ from what is written in ~\cite{Mat2}.

This tree represents the branch locus of the two covers $Rev_1$ and $Rev_2$. The indices used are there to indicate which cover(s) the $"x_i"$ points are branch points of. In particular, the $x_{12}$ points designate branch points common to $Rev_1$ and $Rev_2$. The caps used on the $x_2$ points are just there to distinguish them from each other.

This drawing is only valid for the most general case where $m_1<m_2$ and would be simpler to represent in the case where $m_1=m_2$. 

Figure 2 also shows the other two trees representing the branch locus of each of the two covers $Rev_1$ and $Rev_2$.

 \begin{figure}[h]
\begin{pspicture}(0,-0.5)(16,6)	
%\pscurve[showpoints=true]{-}(0,1.3)
%(0.7,1.8)(3.3,0.5)(4,1.6)(0.4,0.4)
\psline[showpoints=false]{-}(0.4,0.6)(0.45,0.2)(0.65,0.16)(0.25,0.08)(0.65,0.0)(0.25,-0.08)(0.55,-0.14)%generique
\pscurve[showpoints=false]{-}(0,0.25)(0.4,0.4)(4,1.6)(8,1.9)(12,1.6)%horizontale0
\pscurve[showpoints=false]{-}(2,0.5)(1.6,3)(1.6,5)(1.7,6)%verticale1
%\pscurve[showpoints=false]{-}(0.5,4)(2.5,4.8)(4,5.2)(4.5,5.3)%horizontale1
%\pscurve[showpoints=false]{-}(9,1)(8.6,3.5)(8.6,5.5)(8.7,6.5)%verticale2

\psdot(1.73,2)\uput[0](1.73,2){$x_1$}
\psdot(1.58,4.48)\uput[0](1.58,4.48){$x_{12}$}
%\psdot(4,5.2)\uput[90](4,5.2){$x_{1}$}
%\psdot(8.6,3.5)\uput[0](8.6,3.5){$\hat{x}_{2}$}
%\psdot(8.6,5.5)\uput[0](8.6,5.5){$\check{x}_{2}$}

\uput[270](3.5,0.5){$\underbrace{\hspace{4cm}}_{\frac{m_1+1}{2}\;\mathrm{times}}$}
%\uput[270](10.5,0.5){$\underbrace{\hspace{4cm}}_{\frac{m_2-m_1}{2}\;\mathrm{times}}$}

\uput[0](4,3){\Large $\cdots$}%\uput[0](10,4){\Large $\cdots$}
%\psaxes{->}(0,0)(-0.5,-0.5)(14.5,7.5)[$x$,0][$y$, 90]	%creates axes
%\psdot(2,1)	%plots the point (2,1)
%\uput[0](2,1){$A$}	%labels the point (2,1) as A

\end{pspicture}
%\caption{Lieu de branchement  quand $G=(\Z/2\Z)^2$}\label{fig1}

%\end{figure}

%\begin{figure}[h]

\bigskip\bigskip\bigskip
\begin{pspicture}(0,-0.5)(16,6)	
%\pscurve[showpoints=true]{-}(0,1.3)
%(0.7,1.8)(3.3,0.5)(4,1.6)(0.4,0.4)
\psline[showpoints=false]{-}(0.4,0.6)(0.45,0.2)(0.65,0.16)(0.25,0.08)(0.65,0.0)(0.25,-0.08)(0.55,-0.14)%generique
\pscurve[showpoints=false]{-}(0,0.25)(0.4,0.4)(4,1.6)(8,1.9)(12,1.6)%horizontale0
\pscurve[showpoints=false]{-}(2,0.5)(1.6,3)(1.6,5)(1.7,6)%verticale1
%\pscurve[showpoints=false]{-}(0.5,4)(2.5,4.8)(4,5.2)(4.5,5.3)%horizontale1
\pscurve[showpoints=false]{-}(9,1)(8.6,3.5)(8.6,5.5)(8.7,6.5)%verticale2

\psdot(1.73,2)\uput[0](1.73,2){$x_2$}
\psdot(1.58,4.48)\uput[0](1.58,4.48){$x_{12}$}
%\psdot(4,5.2)\uput[90](4,5.2){$x_{1}$}
\psdot(8.6,3.5)\uput[0](8.6,3.5){$\hat{x}_{2}$}
\psdot(8.6,5.5)\uput[0](8.6,5.5){$\check{x}_{2}$}

\uput[270](3.5,0.5){$\underbrace{\hspace{4cm}}_{\frac{m_1+1}{2}\;\mathrm{times}}$}
\uput[270](10.5,0.5){$\underbrace{\hspace{4cm}}_{\frac{m_2-m_1}{2}\;\mathrm{times}}$}

\uput[0](4,3){\Large $\cdots$}\uput[0](10,4){\Large $\cdots$}
%\psaxes{->}(0,0)(-0.5,-0.5)(14.5,7.5)[$x$,0][$y$, 90]	%creates axes
%\psdot(2,1)	%plots the point (2,1)
%\uput[0](2,1){$A$}	%labels the point (2,1) as A

\end{pspicture}
\caption{Branch loci of $Rev_1$ and $Rev_2$}\label{fig2}
$ $
\end{figure}
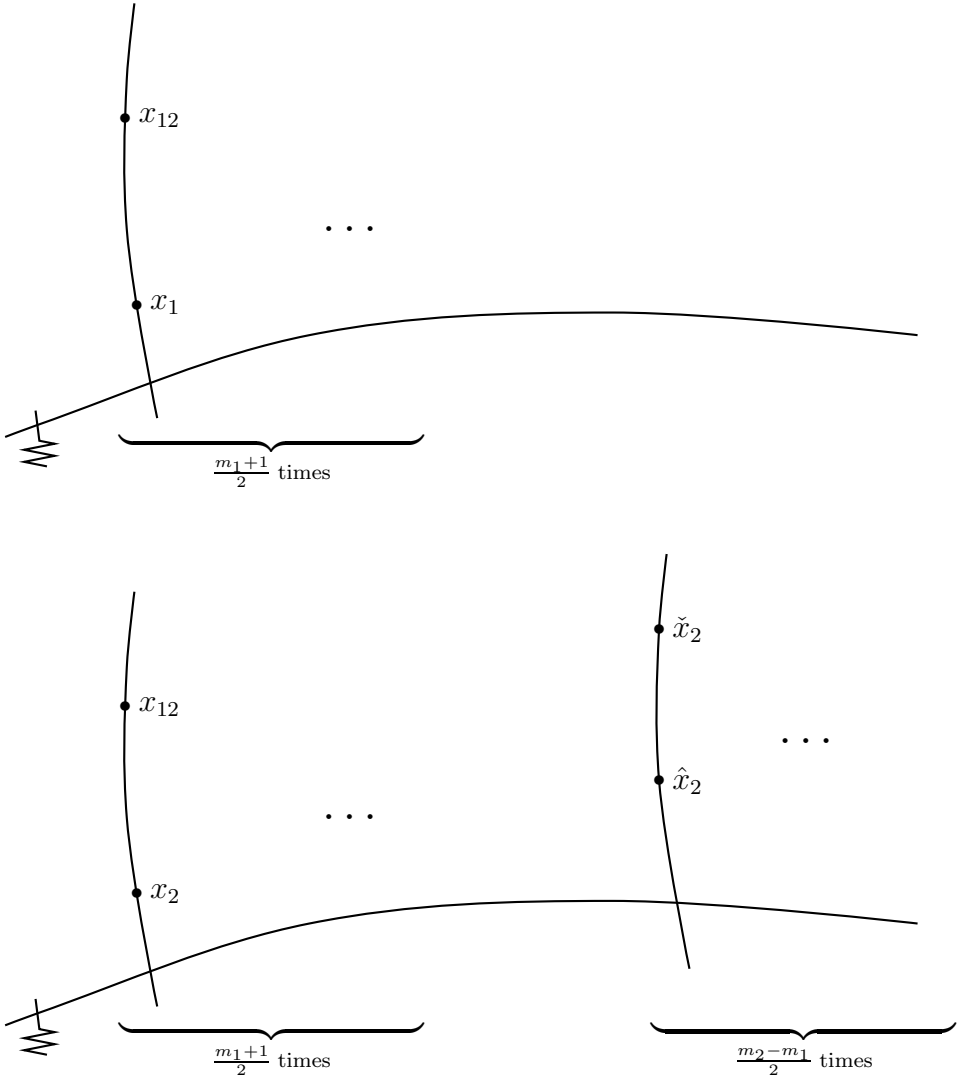

\subsection{Construction of covers with fixed geometry and reduction}
\subsubsection{Aim}

The aim is as follows. Let $r\in\N^*$ and $X_1,\cdots,X_{r}\in R$ be nonzero such that $v(X_i-X_j)=0$ for all $i\neq j$. We consider a $(\Z/2\Z)$-extension of $k[[t]]$ given by an equation of the form

  $$\begin{aligned}
   w^2+w&=\sum_{\ell=0}^{r-1}\frac{a_{\ell}^2}{t^{2r-1-2\ell}}\end{aligned}\;\;\;(E)$$
   where $a_0,\cdots,a_{r-1}\in k$, $a_0\neq 0$. The conductor here is $2r$.

   We then want to construct a polynomial $f\in R[X]$ such that 
  \begin{enumerate}[a.]
      \item The cover of the $2$-adic open-disk given by the equation $Y^2=f(X)$ has good reduction with respect to a certain Gauss valuation, this reduction inducing equation $(E)$. %Pour vérifier cette bonne réduction, on utilisera le théorème 3.4 de ~\cite{GrMa}  : la bonne réduction est obtenue quand on a l'égalité des degrés des différentes générique et spéciale. Pour cela, il suffit de vérifier que le nombre de point de branchement à la fibre générique est égal au conducteur du revêtement d'Artin Schreier à la fibre spéciale.  
    \item The tree describing the geometry of the branch locus has a single internal component and $r$ terminal components, each containing two branch points (this corresponds to the drawing of the first tree in Figure 2, where $\frac{m_1+1}{2}$ is replaced by $r$).
    \item The points $X_1,\cdots,X_{r}$ are part of the branch locus, i.e. $X_1,\cdots,X_{r}$ are simple roots of $f$ and there is exactly one of these points in each terminal component.
  \end{enumerate}

  The idea is to work from the exact differential form on the internal component of the tree, which is of the form 
  $$\frac{dx}{\prod\limits_{i=1}^r(x-x_i)^2}$$
  and then calculate the underlying partial fraction expansion.
\subsubsection{Construction of the cover}

Let's start with a lemma. 

\begin{lemma}
 Let $r\in\N^*$ and $X_1,\cdots,X_{r}\in R$ be non-zero, such that $v(X_i-X_j)=0$ for all $i\neq j$. Let $Q=\prod_{i=1}^r(X-X_i)$. Let $\ell\in\cg0,r-1\cd$. Then there exists $\gamma_{i\ell}\in R$ such that
  $$X^{2\ell}=Q^2(X)\left(\sum_{i=1}^{r}\frac{\gamma_{i\ell}^2}{(X-X_i)^2}\right)+O(2) $$
  where $O(2)$ denotes a polynomial in $2R[X]$ of degree at most $2r-1$.
\end{lemma}

\begin{Proof}\rm
    We write the partial fraction expansion of $\frac{X^{\ell}}{Q}$. There then exists $\gamma_{i\ell}\in K:=Frac(R)$ such that 

     $$\frac{X^{\ell}}{Q(X)}=\sum_{i=1}^{r}\frac{\gamma_{i\ell}}{(X-X_i)}. $$

    The coefficients $\gamma_{i\ell}$ are equal to  $\frac{X_i^{\ell}}{Q'(X_i)}$. They are therefore within $R$ (and even non-zero), since $v(X_i-X_j)=0$ for all $i\neq j$.

     We square this equality and multiply by $Q^2$, which completes the proof. 

     \begin{flushright}
     $\square$
     \end{flushright}
    
\end{Proof}

Let's fix some additional notations. Let $\rho \in R$ be such that $v(2)\leq v(\rho)< 2v(2)$. Let $\rho_0=\left(\frac{4}{\rho}\right)^{\frac{1}{2r-1}}$. Note that $v(\rho_0)>0$.

Let $A_0, \cdots,A_{r-1}\in R$ with $v(A_0)=0$. For all $i\in\cg1,r\cd$ we write :

$$\begin{aligned}
 % R_i(X)&=\sum_{\ell=0}^{n-1}\rho_0^{2\ell}A_{\ell}^2\left(\beta_{i\ell}^2(X-%X_i)^3+\alpha_{i\ell}^2(X-X_i)\right)\\
 % \tilde{R}_i(X)&=\sum_{\ell=0}^{n-1}\rho_0^{\ell}A_{\ell}\left(\beta_{i\ell}(X-X_i)^3+\alpha_{i\ell}(X-X_i)^2\right)X_i^{\frac{1}{2}}\\
  R_i(X)&=\sum_{\ell=0}^{r-1}\rho_0^{2\ell}A_{\ell}^2\gamma_{i\ell}^2(X-X_i)\\
   \tilde{R}_i(X)&=\sum_{\ell=0}^{r-1}\rho_0^{\ell}A_{\ell}\gamma_{i\ell}X_i^{\frac{1}{2}}(X-X_i)\\
   P_i(X)&=(X-X_i)^2+\rho R_i(X)+2\rho^{\frac{1}{2}}\tilde{R}_i(X)
\end{aligned}$$

Finally let $$f(X)=\prod_{i=1}^r P_i(X).$$

We then have the following lemma:
\begin{lemma}
   With the previous notations, we have the relation :
    \begin{align*}
f(X) &= Q^2(X)\left(1+\rho 
\sum_{i=1}^{r}\frac{R_i(X)}{(X-X_i)^2}+2\rho^{\frac{1}{2}}
\sum_{i=1}^{r}\frac{\tilde{R}_i(X)}{(X-X_i)^2}\right)+O(4)
\end{align*}

where $O(4)$ denotes a polynomial in $4R[X]$ of degree at most $2r-1$.
\end{lemma}

\begin{Proof} \rm We have
\begin{align*}  
        f(X)&=Q^2(X)\prod_{i=1}^r \frac{P_i(X)}{(X-X_i)^2}\\
        &=Q^2(X)\prod_{i=1}^r \left[ 1+\rho\frac{R_i(X)}{(X-X_i)^2}+2\rho^{\frac{1}{2}}\frac{\tilde{R}_i(X)}{(X-X_i)^2}\right]\\
        &=Q^2(X) \left[ 1+\rho\sum_{i=1}^r\frac{R_i(X)}{(X-X_i)^2}+2\rho^{\frac{1}{2}}\sum_{i=1}^r\frac{\tilde{R}_i(X)}{(X-X_i)^2}\right]+O(4)
\end{align*}
 since $v(\rho^2)\geq 2v(2)$ and $v(2\rho)\geq 2v(2)$.

   \begin{flushright}
     $\square$
     \end{flushright}
\end{Proof}

We are now in a position to prove the desired result.

\begin{Prop}
Let's consider the same polynomial $f$ as in Lemma 2, and denote $a_0,\cdots,a_{r-1}$ the reductions modulo $\pi$ of $A_0,\cdots,A_{r-1}$. In particular, $a_0\neq 0$. Then the cover of $\mathbb{P}_K^1$ given by the equation $Y^2=f(X)$ has good reduction for the Gauss valuation relative to $T=\rho_0X$, the reduction being :
    $$w^2-w=\sum_{\ell=0}^{r-1}\frac{a_{\ell}^2}{t^{2r-1-2\ell}}.$$

    In addition, the geometry of the branch locus is given by the first tree in figure 2 (with $r$ terminal components instead of $\frac{m_1+1}{2}$).  
   
\end{Prop}

\begin{Proof}\rm We start from the result of Lemma 2. Thus

\begin{align*}  
        f(X)&=Q^2(X) \left[ 1+\rho\sum_{i=1}^r\frac{R_i(X)}{(X-X_i)^2}+2\rho^{\frac{1}{2}}\sum_{i=1}^r\frac{\tilde{R}_i(X)}{(X-X_i)^2}\right]+O(4)\\
        &=Q^2(X) \left[ 1+\rho\sum_{i=1}^r\sum_{\ell=0}^{r-1}\rho_0^{2\ell}A_{\ell}^2\frac{\gamma_{i\ell}^2}{(X-X_i)}+2\rho^{\frac{1}{2}}\sum_{i=1}^r\sum_{\ell=0}^{r-1}\rho_0^{\ell}A_{\ell}\frac{\gamma_{i\ell}}{(X-X_i)}X_i^{\frac{1}{2}}\right]+O(4)\\
        &=Q^2(X) \left[ 1+\rho\sum_{\ell=0}^{r-1}\sum_{i=1}^r\rho_0^{2\ell}A_{\ell}^2\frac{\gamma_{i\ell}^2(X-X_i)}{(X-X_i)^2}+2\rho^{\frac{1}{2}}\sum_{\ell=0}^{r-1}\sum_{i=1}^r\rho_0^{\ell}A_{\ell}\frac{\gamma_{i\ell}}{(X-X_i)}X_i^{\frac{1}{2}}\right]+O(4)\\
        &=Q^2(X)+\rho XQ^2(X)\sum_{\ell=0}^{r-1}\rho_0^{2\ell}A_{\ell}^2\sum_{i=1}^r\frac{\gamma_{i\ell}^2}{(X-X_i)^2}-\rho Q^2(X)\sum_{\ell=0}^{r-1}\rho_0^{2\ell}A_{\ell}^2\sum_{i=1}^r\frac{\gamma_{i\ell}^2X_i}{(X-X_i)^2}\\
        &+2\rho^{\frac{1}{2}}Q^2(X)\sum_{\ell=0}^{r-1}\sum_{i=1}^r\rho_0^{\ell}A_{\ell}\frac{\gamma_{i\ell}}{(X-X_i)}X_i^{\frac{1}{2}}+O(4)\\
        &=Q^2(X)\left[1+\rho\left(\sum_{\ell=0}^{r-1}\rho_0^{\ell}A_{\ell}\sum_{i=1}^r\frac{\gamma_{i\ell}X_i^{\frac{1}{2}}}{(X-X_i)}\right)^2+2\rho^{\frac{1}{2}}\sum_{\ell=0}^{r-1}\sum_{i=1}^r\rho_0^{\ell}A_{\ell}\frac{\gamma_{i\ell}}{(X-X_i)}X_i^{\frac{1}{2}}\right]\\
        &+\rho X\sum_{\ell=0}^{r-1}\rho_0^{2\ell}A_{\ell}^2Q^2(X)\sum_{i=1}^r\frac{\gamma_{i\ell}^2}{(X-X_i)^2}+O(4)
\end{align*}

since $v(2\rho)\geq 2v(2)$. Now let's use Lemma 1. We get :
\begin{align*}
    f(X)&=Q^2(X)\left[1+\rho^{\frac{1}{2}}\sum_{\ell=0}^{r-1}\rho_0^{\ell}A_{\ell}\frac{\gamma_{i\ell}}{(X-X_i)}X_i^{\frac{1}{2}}\right]^2+\rho\sum_{\ell=0}^{r-1}\rho_0^{2\ell}A_{\ell}^2X^{2\ell+1}+O(4)\\
    &=Q_0^2(X)+\rho\sum_{\ell=0}^{r-1}\rho_0^{2\ell}A_{\ell}^2X^{2\ell+1}+O(4)
\end{align*}
    assuming $Q_0(X)=Q(X)\left[1+\rho^{\frac{1}{2}}\sum_{\ell=0}^{r-1}\rho_0^{\ell}A_{\ell}\frac{\gamma_{i\ell}}{(X-X_i)}X_i^{\frac{1}{2}}\right]$.

  Let $Y=2V+Q_0(X)$. The equation $Y^2=f(X)$ then becomes
  $$4(V^2+Q_0(X)V)= \rho\sum_{\ell=0}^{r-1}\rho_0^{2\ell}A_{\ell}^2X^{2\ell+1}+O(4) $$ i.e.
  $$ \frac{4(V^2+Q_0(X)V)}{X^{2r}}=\rho\rho_0^{2r-1}\sum_{\ell=0}^{r-1}A_{\ell}^2\left(\frac{1}{\rho_0X}\right)^{2r-1-2\ell}+\frac{O(4)}{X^{2r}}.$$

  But $\rho\rho_0^{2r-1}=4$, then we get
  $$ \frac{(V^2+Q_0(X)V)}{X^{2r}}=\sum_{\ell=0}^{r-1}A_{\ell}^2\left(\frac{1}{\rho_0X}\right)^{2r-1-2\ell}+\frac{O(1)}{X^{2r}}$$

  Here $O(1)$ denotes an element of $R[X]$ of degree at most $2r-1$. 
  Let's write $T=\rho_0X$ and $V=X^rW$. The equation becomes
  $$W^2+W\frac{Q_0(X)}{X^{r}}=\sum_{\ell=0}^{r-1}A_{\ell}^2\left(\frac{1}{T}\right)^{2r-1-2\ell}+\frac{O(1)}{X^{2r}}, $$
  which gives in reduction 
  $$w^2-w=\sum_{\ell=0}^{r-1}\frac{a_{\ell}^2}{t^{2r-1-2\ell}}.$$

It's easy to check that the geometry of the branch locus is as claimed., since $X_i$ is one of the roots of $P_i$ and the other root (let's call it $X'_i$) is such that $v(X_i-X'_i)=v(\rho)$.
   \begin{flushright}
     $\square$
     \end{flushright}
\end{Proof}

\subsection{End of the proof of Theorem 1.a}

Consider a $(\Z/2\Z)^2$-extension of $k[[t]]$ given by the equations 
$$\left\{
\begin{aligned}
    w_1^2-w_1&=\sum_{\ell=0}^{\frac{m_1-1}{2}}\frac{a_{\ell}^2}{t^{m_1-2\ell}}\;\;\;\;(1)\\
    w_2^2-w_2&=\sum_{\ell=0}^{\frac{m_2-1}{2}}\frac{b_{\ell}^2}{t^{m_2-2\ell}}\;\;\;\;(2)
\end{aligned}
\right.$$
where $a_{\ell}$, $b_{\ell} \in k$ and $a_0,b_0\in k^*$.

Let $X_1,\cdots, X_{\frac{m_2+1}{2}}\in R$, non-zero, such that $v(X_i-X_j)=0$ if $i\neq j$. We'll apply Proposition 1 twice, with judicious choices for $r$ and $\rho$.
\begin{enumerate}[$\bullet$]
    \item Let's first apply Proposition 1 to the case where $r=\frac{m_1+1}{2}$ and $\rho=2^{2-\frac{m_1}{m_2}}$ and the $A_{\ell}$ lifting the $a_{\ell}$. For this, it's easy to check that the condition $v(2)\leq v(\rho)<2v(2)$ is satisfied. We then find a cover of equation $Y^2=f_1(X)$ which has good reduction for the $T$-Gauss valuation (where $T=\rho_0X$). 
    In this case, we have 
    $$ \rho_0=\left(\frac{4}{\rho}\right)^{\frac{1}{2r-1}}=
    \left(2^{\frac{m_1}{m_2}}\right)^{\frac{1}{m_1}}=2^{\frac{1}{m_2}}.$$
    and the reduction gives the equation (1).
    \item Let's apply Proposition 1 a second time in the case where $r=\frac{m_2+1}{2}$ and $\rho=2$ and the $A_{\ell}$ lifting the $b_{\ell}$. For this, it's easy to check that the condition $v(2)\leq v(\rho)<2v(2)$ is satisfied. We then find a cover of equation $Y^2=f_2(X)$ which has good reduction  for the $T$-Gauss valuation (where $T=\rho_0X$). 
    This time we have 
    $$ \rho_0=\left(\frac{4}{\rho}\right)^{\frac{1}{2r-1}}=2^{\frac{1}{m_2}}$$
    and the reduction gives the equation (2).
\end{enumerate}

The two covers therefore have simultaneously good reduction and have at least $\frac{m_1+1}{2}$ branch points in common (these are the $X_1,\cdots, X_{\frac{m_1+1}{2}}$). 

In the case $m_1<m_2$, we can be sure that there can be no more points in common, as the radii of the disks corresponding to the terminal components of the branch tree of $Rev_1$ are not the same as those of $Rev_2$.

In the case where $m_1=m_2$, Kato's result (see assertion 3.4 in ~\cite{GrMa}) then ensures that there can be no more branch points in common (the degree of the generic different is always greater than that of the special different). 

In any case we have exactly $\frac{m_1+1}{2}$ branch points in common, which using Theorem 2 completes the proof.
   
\section{The case $(\Z/2\Z)^3$ }
  \subsection{Geometry of the branch locus}
   First, we want to generalize the geometry of the branch locus as constructed in Section IV to the case $G=(\Z/2\Z)^3$.
    We therefore start with an extension $k[[z]]/\!\raisebox{-.65ex}{\ensuremath{k[[z]]^G}}$, of type $(m_1+1,m_2+1,m_3+1)$. Recall that the only condition on $m_i$ is that $4|m_1+1$. Moreover, according to Theorem 2, the combinatorics of the branch loci of the three covers $Rev_1$, $Rev_2$, $Rev_3$ must be such that :
   \begin{enumerate}[-]
       \item $Rev_1$, $Rev_2$, $Rev_3$ have $\frac{m_1+1}{4}$ branch points in common.
       \item$Rev_1$ and $Rev_2$ have $\frac{m_1+1}{2}$ branch points in common. 
       \item $Rev_1$ and $Rev_3$ have $\frac{m_1+1}{2}$ branch points in common.
       \item $Rev_2$ and $Rev_3$ have $\frac{m_2+1}{2}$ branch points in common.
    \end{enumerate}

 The idea is to consider the geometry given in figure 3.
  $ $
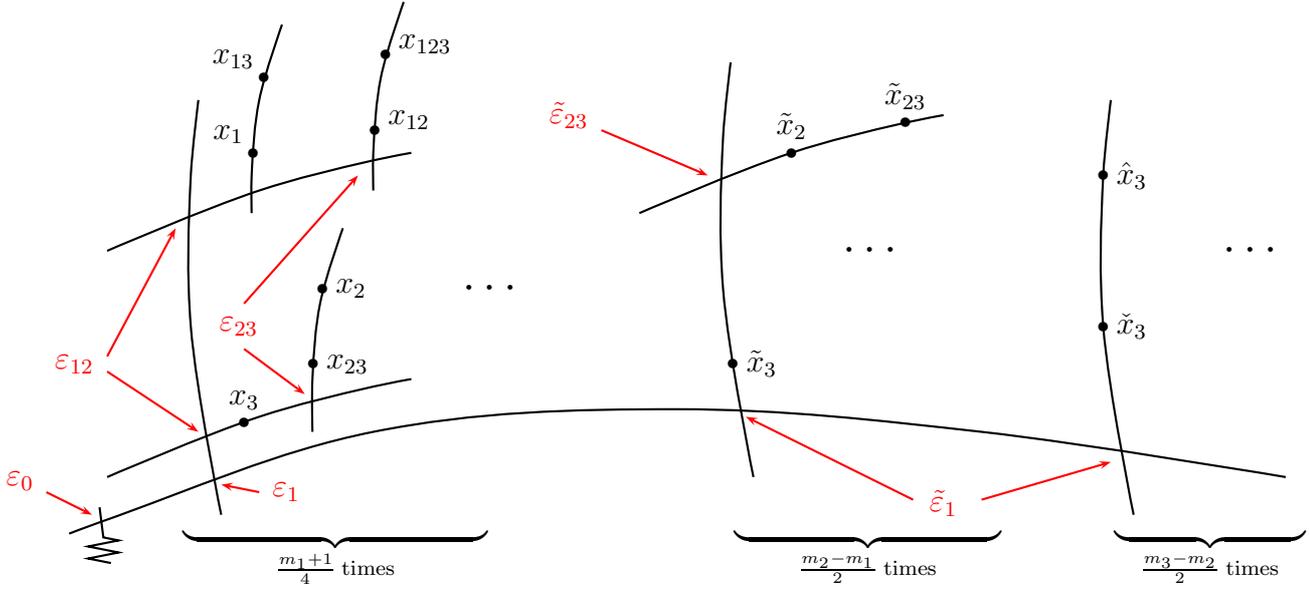
\begin{figure}[h]
\begin{pspicture}(0,-0.5)(18,8.5)	
%\pscurve[showpoints=true]{-}(0,1.3)
%(0.7,1.8)(3.3,0.5)(4,1.6)(0.4,0.4)
\psline[showpoints=false]{-}(0.4,0.6)(0.45,0.2)(0.65,0.16)(0.25,0.08)(0.65,0.0)(0.25,-0.08)(0.55,-0.14)%generique
\pscurve[showpoints=false]{-}(0,0.25)(0.4,0.4)(4,1.6)(8,1.9)(12,1.6)(16,1)%horizontale0
\pscurve[showpoints=false]{-}(2,0.5)(1.6,3)(1.6,5)(1.7,6)%verticale1
\pscurve[showpoints=false]{-}(2.4,4.5)(2.4,5)(2.5,6)(2.8,7)%verticale1bis
\pscurve[showpoints=false]{-}(4,4.8)(4,5.3)(4.1,6.3)(4.4,7.3)%verticale1ter
\pscurve[showpoints=false]{-}(3.2,1.6)(3.2,2.3)(3.3,3.3)(3.6,4.3)%verticale4
\pscurve[showpoints=false]{-}(0.5,4)(2.5,4.8)(4,5.2)(4.5,5.3)%horizontale1
\pscurve[showpoints=false]{-}(0.5,1)(2.5,1.8)(4,2.2)(4.5,2.3)%horizontale2
\pscurve[showpoints=false]{-}(9,1)(8.6,3.5)(8.6,5.5)(8.7,6.5)%verticale2

\pscurve[showpoints=false]{-}(7.5,4.5)(9.5,5.3)(11,5.7)(11.5,5.8)%horizontale1
\psdot(8.73,2.5)\uput[0](8.73,2.5){$\tilde{x}_3$}
\psdot(9.5,5.3)\uput[90](9.5,5.3){$\tilde{x}_{2}$}
\psdot(11,5.7)\uput[90](11,5.7){$\tilde{x}_{23}$}

\psdot(2.3,1.72)\uput[90](2.3,1.72){$x_3$}
\psdot(3.21,2.5)\uput[0](3.21,2.5){$x_{23}$}
\psdot(3.33,3.5)\uput[0](3.33,3.5){$x_{2}$}
\psdot(2.42,5.3)\uput[135](2.42,5.3){$x_1$}
\psdot(2.56,6.3)\uput[135](2.56,6.3){$x_{13}$}
\psdot(4.02,5.6)\uput[25](4.02,5.6){$x_{12}$}
\psdot(4.16,6.6)\uput[25](4.16,6.6){$x_{123}$}

\pscurve[showpoints=false]{-}(14,0.5)(13.6,3)(13.6,5)(13.7,6)%verticale1
\psdot(13.6,3)\uput[0](13.6,3){$\check{x}_{3}$}
\psdot(13.6,5)\uput[0](13.6,5){$\hat{x}_{3}$}
%\psdot(8.6,3.5)\uput[0](8.6,3.5){$\hat{x}_{2}$}
%\psdot(8.6,5.5)\uput[0](8.6,5.5){$\check{x}_{2}$}

\uput[270](3.5,0.5){$\underbrace{\hspace{4cm}}_{\frac{m_1+1}{4}\;\mathrm{times}}$}
\uput[270](10.5,0.5){$\underbrace{\hspace{3.5cm}}_{\frac{m_2-m_1}{2}\;\mathrm{times}}$}
\uput[270](15,0.5){$\underbrace{\hspace{2.5cm}}_{\frac{m_3-m_2}{2}\;\mathrm{times}}$}

\uput[0](5,3.5){\Large $\cdots$}\uput[0](10,4){\Large $\cdots$}\uput[0](15,4){\Large $\cdots$}
%\psaxes{->}(0,0)(-0.5,-0.5)(14.5,7.5)[$x$,0][$y$, 90]	%creates axes
%\psdot(2,1)	%plots the point (2,1)
%\uput[0](2,1){$A$}	%labels the point (2,1) as A

\uput[350](-1,1){\textcolor{red}{$\varepsilon_0$}}\psline[linecolor=red]{->}(-0.3,0.8)(0.3,0.5)
\uput[0](2.5,0.8){\textcolor{red}{$\varepsilon_1$}}\psline[linecolor=red]{->}(2.5,0.8)(2,0.9)
\uput[180](0.5,2.5){\textcolor{red}{$\varepsilon_{12}$}}\psline[linecolor=red]{->}(0.5,2.4)(1.7,1.6)\psline[linecolor=red]{->}(0.5,2.6)(1.4,4.3)
\uput[0](1.8,3){\textcolor{red}{$\varepsilon_{23}$}}\psline[linecolor=red]{->}(2.3,2.7)(3.1,2.1)\psline[linecolor=red]{->}(2.3,3.3)(3.8,5)
\uput[270](11.5,1){\textcolor{red}{$\tilde{\varepsilon}_1$}}\psline[linecolor=red]{->}(11.1,0.7)(8.9,1.8)\psline[linecolor=red]{->}(12,0.7)(13.7,1.2)
\uput[180](7,5.8){\textcolor{red}{$\tilde{\varepsilon}_{23}$}}
\psline[linecolor=red]{->}(7,5.6)(8.4,5)
\end{pspicture}
\caption{Branch locus when $G=(\Z/2\Z)^3$}\label{fig3}
\end{figure}

  The figure represents the most general case $m_1<m_2<m_3$. Once again, the indices of the points indicate which covers these points are the branch points of. In this way, we can verify that the combinatorics of branch points is respected. In red, we've also given a name to the thicknesses of the different annuli corresponding to the double points. These thicknesses need to be "well chosen" and this is the subject of the next section.

  \subsection{Choice of thicknesses}

The thicknesses must be chosen according to the constraints linked to the variation of the different. In other words, for each branch point, equality (**) of Section III must be verified.  
%Formulated differently, this means using the metric properties of Hurwitz %trees attached to a $\sigma$-automorphism of order $p$ of $R[[Z]]$. More %precisely, the branching tree is the quotient of the Hurwitz tree attached %to an automorphism $\sigma$ of order $p$. For further details, see %Proposition 1.1 in ~\cite{GrMa2}.

      Starting from the point $x_{123}$ and considering successively the three covers $Rev_1$, $Rev_2$, $Rev_3$, we must have :

      $$\left\{
       \begin{aligned}
           \varepsilon_{23}+3(\varepsilon_{12}+\varepsilon_{1})+m_1\varepsilon_{0}&=2v(2)\\
           \varepsilon_{23}+3\varepsilon_{1}+m_2\varepsilon_{0}&=2v(2)\\
           \varepsilon_{12}+3\varepsilon_{1}+m_3\varepsilon_{0}&=2v(2)
       \end{aligned}
      \right.
      $$

      The four unknown thicknesses satisfy three equations, giving us one degree of freedom. Let's arbitrarily choose $\varepsilon_1=\frac{1}{2}v(2)$. Solving the system gives

    $$ \varepsilon_{0}=\frac{1}{m}v(2),\;\;\;\varepsilon_{12}=\frac{m_2-m_1}{2m}v(2),\;\;\;\varepsilon_{23}=\frac{m_3-m_2}{m}v(2)$$
where $m=2m_3+m_2-m_1$.

We reason in the same way to find the thicknesses $\tilde{\varepsilon}_{1}$ and $\tilde{\varepsilon}_{23}$.We find
 $$ \Tilde{\varepsilon}_{1}=\frac{3m_3+2m_2-2m_1}{m}v(2),\;\;\;\tilde{\varepsilon}_{23}=\varepsilon_{23}.$$

Let
$$\varrho_1:=2^{\frac{\varepsilon_1}{v(2)}}=2^{\frac{1}{2}},\;\;
\varrho_2:=2^{\frac{\varepsilon_1+\varepsilon_{12}}{v(2)}}=2^{\frac{m_3+m_2-m_1}{m}},\;\;
\varrho_3:=2^{\frac{\varepsilon_1+\varepsilon_{12}+\varepsilon_{23}}{v(2)}}=2^{\frac{2m_3-m_1}{m}}.$$

It is then clear that $v(\varrho_1)\leq v(\varrho_2)\leq v(\varrho_3)$. These elements then appear in Figure 4, where we've indicated the branch loci of each of the three covers $Rev_1$, $Rev_2$, $Rev_3$. We have also indicated the radius values of the various disks corresponding to the tree components.

\begin{figure}[H]
\begin{pspicture}(0,-0.5)(18,6)	
%\pscurve[showpoints=true]{-}(0,1.3)
%(0.7,1.8)(3.3,0.5)(4,1.6)(0.4,0.4)
\psline[showpoints=false]{-}(0.4,0.6)(0.45,0.2)(0.65,0.16)(0.25,0.08)(0.65,0.0)(0.25,-0.08)(0.55,-0.14)%generique
\pscurve[showpoints=false]{-}(0,0.25)(0.4,0.4)(4,1.6)(8,1.9)(12,1.6)(16,1)%horizontale0
\pscurve[showpoints=false]{-}(2,0.5)(1.6,3)(1.6,5)(1.7,6)%verticale1
%\pscurve[showpoints=false]{-}(2.4,2)(2.4,2.5)(2.5,3.5)(2.8,4.5)%verticale1bis
%\pscurve[showpoints=false]{-}(4,2.3)(4,2.8)(4.1,3.8)(4.4,4.8)%verticale1ter
%\pscurve[showpoints=false]{-}(3.2,1.6)(3.2,2.3)(3.3,3.3)(3.6,4.3)%verticale4
\pscurve[showpoints=false]{-}(0.5,1.5)(2.5,2.3)(4,2.7)(4.5,2.8)%horizontale1
\pscurve[showpoints=false]{-}(0.5,4)(2.5,4.8)(4,5.2)(4.5,5.3)%horizontale2
%\pscurve[showpoints=false]{-}(0.5,1)(2.5,1.8)(4,2.2)(4.5,2.3)%horizontale2
%\pscurve[showpoints=false]{-}(9,1)(8.6,3.5)(8.6,5.5)(8.7,6.5)%verticale2

%\pscurve[showpoints=false]{-}(7.5,4.5)(9.5,5.3)(11,5.7)(11.5,5.8)%horizontale1
%\psdot(8.73,2.5)\uput[0](8.73,2.5){$\tilde{x}_3$}
%\psdot(9.5,5.3)\uput[90](9.5,5.3){$\tilde{x}_{2}$}
%\psdot(11,5.7)\uput[90](11,5.7){$\tilde{x}_{23}$}

%\psdot(2.3,1.72)\uput[90](2.3,1.72){$x_3$}
%\psdot(3.21,2.5)\uput[0](3.21,2.5){$x_{23}$}
%\psdot(3.33,3.5)\uput[0](3.33,3.5){$x_{2}$}
\psdot(2.5,2.3)\uput[135](2.5,2.3){$x_{12}$}
\psdot(4,2.7)\uput[135](4,2.7){$x_{123}$}
\psdot(2.5,4.8)\uput[135](2.5,4.8){$x_{1}$}
\psdot(4,5.2)\uput[135](4,5.2){$x_{13}$}

\uput[0](5,4){\textcolor{red}{$v(\varrho_3)$}}\psline[linecolor=red]{->}(5.3,4.3)(4.8,5.1)\psline[linecolor=red]{->}(5.3,3.6)(4.8,2.8)
\uput[180](0,3){\textcolor{red}{$v(\varrho_2)$}}\psline[linecolor=red]{->}(0,3)(1.2,3)
\uput[90](16,3){\textcolor{red}{$0$}}\psline[linecolor=red]{->}(16,3)(15.5,1.3)

%\pscurve[showpoints=false]{-}(14,0.5)(13.6,3)(13.6,5)(13.7,6)%verticale1
%\psdot(13.6,3)\uput[0](13.6,3){$\check{x}_{3}$}
%\psdot(13.6,5)\uput[0](13.6,5){$\hat{x}_{3}$}
%\psdot(8.6,3.5)\uput[0](8.6,3.5){$\hat{x}_{2}$}
%\psdot(8.6,5.5)\uput[0](8.6,5.5){$\check{x}_{2}$}

\uput[270](3.5,0.5){$\underbrace{\hspace{4cm}}_{\frac{m_1+1}{4}\;\mathrm{times}}$}
%\uput[270](10.5,0.5){$\underbrace{\hspace{3.5cm}}_{\frac{m_2-m_1}{2}\;\mathrm{times}}$}
%\uput[270](15,0.5){$\underbrace{\hspace{2.5cm}}_{\frac{m_3-m_2}{2}\;\mathrm{times}}$}

%\uput[0](7,3.5){\Large $\cdots$}%\uput[0](10,4){\Large $\cdots$}\uput[0](15,4){\Large $\cdots$}
%\psaxes{->}(0,0)(-0.5,-0.5)(14.5,7.5)[$x$,0][$y$, 90]	%creates axes
%\psdot(2,1)	%plots the point (2,1)
%\uput[0](2,1){$A$}	%labels the point (2,1) as A

%\uput[350](-1,1){\textcolor{red}{$\varepsilon_0$}}\psline[linecolor=red]{->}(-0.3,0.8)(0.3,0.5)
%\uput[0](2.5,0.8){\textcolor{red}{$\varepsilon_1$}}\psline[linecolor=red]{->}(2.5,0.8)(2,0.9)
%\uput[180](0.5,2.5){\textcolor{red}{$\varepsilon_{12}$}}\psline[linecolor=red]{->}(0.5,2.4)(1.7,1.6)\psline[linecolor=red]{->}(0.5,2.6)(1.4,4.3)
%\uput[0](1.8,3){\textcolor{red}{$\varepsilon_{23}$}}\psline[linecolor=red]{->}(2.3,2.7)(3.1,2.1)\psline[linecolor=red]{->}(2.3,3.3)(3.8,5)
%\uput[270](11.5,1){\textcolor{red}{$\tilde{\varepsilon}_1$}}\psline[linecolor=red]{->}(11.1,0.7)(8.9,1.8)\psline[linecolor=red]{->}(12,0.7)(13.7,1.2)
%\uput[180](7,5.8){\textcolor{red}{$\tilde{\varepsilon}_{23}$}}\psline[linecolor=red]{->}(7,5.6)(8.4,5)
\end{pspicture}
$ $
\begin{pspicture}(0,-0.5)(18,8)	
%\pscurve[showpoints=true]{-}(0,1.3)
%(0.7,1.8)(3.3,0.5)(4,1.6)(0.4,0.4)
\psline[showpoints=false]{-}(0.4,0.6)(0.45,0.2)(0.65,0.16)(0.25,0.08)(0.65,0.0)(0.25,-0.08)(0.55,-0.14)%generique
\pscurve[showpoints=false]{-}(0,0.25)(0.4,0.4)(4,1.6)(8,1.9)(12,1.6)(16,1)%horizontale0
\pscurve[showpoints=false]{-}(2,0.5)(1.6,3)(1.6,5)(1.7,6)%verticale1
%\pscurve[showpoints=false]{-}(2.4,2)(2.4,2.5)(2.5,3.5)(2.8,4.5)%verticale1bis
%\pscurve[showpoints=false]{-}(4,2.3)(4,2.8)(4.1,3.8)(4.4,4.8)%verticale1ter
%\pscurve[showpoints=false]{-}(3.2,1.6)(3.2,2.3)(3.3,3.3)(3.6,4.3)%verticale4
\pscurve[showpoints=false]{-}(0.5,1.5)(2.5,2.3)(4,2.7)(4.5,2.8)%horizontale1
\pscurve[showpoints=false]{-}(0.5,4)(2.5,4.8)(4,5.2)(4.5,5.3)%horizontale2
%\pscurve[showpoints=false]{-}(0.5,1)(2.5,1.8)(4,2.2)(4.5,2.3)%horizontale2
\pscurve[showpoints=false]{-}(9,1)(8.6,3.5)(8.6,5.5)(8.7,6.5)%verticale2

\psdot(8.6,3.5)\uput[0](8.6,3.5){$\tilde{x}_2$}
\psdot(8.6,5.5)\uput[0](8.6,5.5){$\tilde{x}_{23}$}
\uput[90](11,4.5){\textcolor{red}{$v(\varrho_1^2\varrho_3)$}}\psline[linecolor=red]{->}(10.3,4.8)(8.9,4.3)

%\pscurve[showpoints=false]{-}(7.5,4.5)(9.5,5.3)(11,5.7)(11.5,5.8)%horizontale1
%\psdot(8.73,2.5)\uput[0](8.73,2.5){$\tilde{x}_3$}
%\psdot(9.5,5.3)\uput[90](9.5,5.3){$\tilde{x}_{2}$}
%\psdot(11,5.7)\uput[90](11,5.7){$\tilde{x}_{23}$}

%\psdot(2.3,1.72)\uput[90](2.3,1.72){$x_3$}
%\psdot(3.21,2.5)\uput[0](3.21,2.5){$x_{23}$}
%\psdot(3.33,3.5)\uput[0](3.33,3.5){$x_{2}$}
\psdot(2.5,2.3)\uput[135](2.5,2.3){$x_{12}$}
\psdot(4,2.7)\uput[135](4,2.7){$x_{123}$}
\psdot(2.5,4.8)\uput[135](2.5,4.8){$x_{2}$}
\psdot(4,5.2)\uput[135](4,5.2){$x_{23}$}

\uput[0](5,4){\textcolor{red}{$v(\varrho_3)$}}\psline[linecolor=red]{->}(5.3,4.3)(4.8,5.1)\psline[linecolor=red]{->}(5.3,3.6)(4.8,2.8)
\uput[180](0,3){\textcolor{red}{$v(\varrho_1)$}}\psline[linecolor=red]{->}(0,3)(1.2,3)
\uput[90](16,3){\textcolor{red}{$0$}}\psline[linecolor=red]{->}(16,3)(15.5,1.3)

%\pscurve[showpoints=false]{-}(14,0.5)(13.6,3)(13.6,5)(13.7,6)%verticale1
%\psdot(13.6,3)\uput[0](13.6,3){$\check{x}_{3}$}
%\psdot(13.6,5)\uput[0](13.6,5){$\hat{x}_{3}$}
%\psdot(8.6,3.5)\uput[0](8.6,3.5){$\hat{x}_{2}$}
%\psdot(8.6,5.5)\uput[0](8.6,5.5){$\check{x}_{2}$}

\uput[270](3.5,0.5){$\underbrace{\hspace{4cm}}_{\frac{m_1+1}{4}\;\mathrm{times}}$}
\uput[270](10.5,0.5){$\underbrace{\hspace{3.5cm}}_{\frac{m_2-m_1}{2}\;\mathrm{times}}$}
%\uput[270](15,0.5){$\underbrace{\hspace{2.5cm}}_{\frac{m_3-m_2}{2}\;\mathrm{times}}$}

%\uput[0](7,3.5){\Large $\cdots$}%\uput[0](10,4){\Large $\cdots$}\uput[0](15,4){\Large $\cdots$}
%\psaxes{->}(0,0)(-0.5,-0.5)(14.5,7.5)[$x$,0][$y$, 90]	%creates axes
%\psdot(2,1)	%plots the point (2,1)
%\uput[0](2,1){$A$}	%labels the point (2,1) as A

%\uput[350](-1,1){\textcolor{red}{$\varepsilon_0$}}\psline[linecolor=red]{->}(-0.3,0.8)(0.3,0.5)
%\uput[0](2.5,0.8){\textcolor{red}{$\varepsilon_1$}}\psline[linecolor=red]{->}(2.5,0.8)(2,0.9)
%\uput[180](0.5,2.5){\textcolor{red}{$\varepsilon_{12}$}}\psline[linecolor=red]{->}(0.5,2.4)(1.7,1.6)\psline[linecolor=red]{->}(0.5,2.6)(1.4,4.3)
%\uput[0](1.8,3){\textcolor{red}{$\varepsilon_{23}$}}\psline[linecolor=red]{->}(2.3,2.7)(3.1,2.1)\psline[linecolor=red]{->}(2.3,3.3)(3.8,5)
%\uput[270](11.5,1){\textcolor{red}{$\tilde{\varepsilon}_1$}}\psline[linecolor=red]{->}(11.1,0.7)(8.9,1.8)\psline[linecolor=red]{->}(12,0.7)(13.7,1.2)
%\uput[180](7,5.8){\textcolor{red}{$\tilde{\varepsilon}_{23}$}}\psline[linecolor=red]{->}(7,5.6)(8.4,5)
\end{pspicture}
\begin{pspicture}(0,-0.5)(18,8)	
%\pscurve[showpoints=true]{-}(0,1.3)
%(0.7,1.8)(3.3,0.5)(4,1.6)(0.4,0.4)
\psline[showpoints=false]{-}(0.4,0.6)(0.45,0.2)(0.65,0.16)(0.25,0.08)(0.65,0.0)(0.25,-0.08)(0.55,-0.14)%generique
\pscurve[showpoints=false]{-}(0,0.25)(0.4,0.4)(4,1.6)(8,1.9)(12,1.6)(16,1)%horizontale0
\pscurve[showpoints=false]{-}(2,0.5)(1.6,3)(1.6,5)(1.7,6)%verticale1
%\pscurve[showpoints=false]{-}(2.4,2)(2.4,2.5)(2.5,3.5)(2.8,4.5)%verticale1bis
%\pscurve[showpoints=false]{-}(4,2.3)(4,2.8)(4.1,3.8)(4.4,4.8)%verticale1ter
%\pscurve[showpoints=false]{-}(3.2,1.6)(3.2,2.3)(3.3,3.3)(3.6,4.3)%verticale4
\pscurve[showpoints=false]{-}(0.5,1.5)(2.5,2.3)(4,2.7)(4.5,2.8)%horizontale1
\pscurve[showpoints=false]{-}(0.5,4)(2.5,4.8)(4,5.2)(4.5,5.3)%horizontale2
%\pscurve[showpoints=false]{-}(0.5,1)(2.5,1.8)(4,2.2)(4.5,2.3)%horizontale2
\pscurve[showpoints=false]{-}(9,1)(8.6,3.5)(8.6,5.5)(8.7,6.5)%verticale2

\psdot(8.6,3.5)\uput[0](8.6,3.5){$\tilde{x}_3$}
\psdot(8.6,5.5)\uput[0](8.6,5.5){$\tilde{x}_{23}$}
\uput[90](11,4.5){\textcolor{red}{$v(\varrho_1^2\varrho_2)$}}\psline[linecolor=red]{->}(10.3,4.8)(8.9,4.3)\psline[linecolor=red]{->}(11.8,4.8)(13.3,4.3)

%\pscurve[showpoints=false]{-}(7.5,4.5)(9.5,5.3)(11,5.7)(11.5,5.8)%horizontale1
%\psdot(8.73,2.5)\uput[0](8.73,2.5){$\tilde{x}_3$}
%\psdot(9.5,5.3)\uput[90](9.5,5.3){$\tilde{x}_{2}$}
%\psdot(11,5.7)\uput[90](11,5.7){$\tilde{x}_{23}$}

%\psdot(2.3,1.72)\uput[90](2.3,1.72){$x_3$}
%\psdot(3.21,2.5)\uput[0](3.21,2.5){$x_{23}$}
%\psdot(3.33,3.5)\uput[0](3.33,3.5){$x_{2}$}
\psdot(2.5,2.3)\uput[135](2.5,2.3){$x_{13}$}
\psdot(4,2.7)\uput[135](4,2.7){$x_{123}$}
\psdot(2.5,4.8)\uput[135](2.5,4.8){$x_{3}$}
\psdot(4,5.2)\uput[135](4,5.2){$x_{23}$}

\uput[0](5,4){\textcolor{red}{$v(\varrho_2)$}}\psline[linecolor=red]{->}(5.3,4.3)(4.8,5.1)\psline[linecolor=red]{->}(5.3,3.6)(4.8,2.8)
\uput[180](0,3){\textcolor{red}{$v(\varrho_1)$}}\psline[linecolor=red]{->}(0,3)(1.2,3)
\uput[90](16,3){\textcolor{red}{$0$}}\psline[linecolor=red]{->}(16,3)(15.5,1.3)

\pscurve[showpoints=false]{-}(14,0.5)(13.6,3)(13.6,5)(13.7,6)%verticale1
\psdot(13.6,3)\uput[0](13.6,3){$\check{x}_{3}$}
\psdot(13.6,5)\uput[0](13.6,5){$\hat{x}_{3}$}
%\psdot(8.6,3.5)\uput[0](8.6,3.5){$\hat{x}_{2}$}
%\psdot(8.6,5.5)\uput[0](8.6,5.5){$\check{x}_{2}$}

\uput[270](3.5,0.5){$\underbrace{\hspace{4cm}}_{\frac{m_1+1}{4}\;\mathrm{times}}$}
\uput[270](10.5,0.5){$\underbrace{\hspace{3.5cm}}_{\frac{m_2-m_1}{2}\;\mathrm{times}}$}
\uput[270](15,0.5){$\underbrace{\hspace{2.5cm}}_{\frac{m_3-m_2}{2}\;\mathrm{times}}$}

%\uput[0](7,3.5){\Large $\cdots$}%\uput[0](10,4){\Large $\cdots$}\uput[0](15,4){\Large $\cdots$}
%\psaxes{->}(0,0)(-0.5,-0.5)(14.5,7.5)[$x$,0][$y$, 90]	%creates axes
%\psdot(2,1)	%plots the point (2,1)
%\uput[0](2,1){$A$}	%labels the point (2,1) as A

%\uput[350](-1,1){\textcolor{red}{$\varepsilon_0$}}\psline[linecolor=red]{->}(-0.3,0.8)(0.3,0.5)
%\uput[0](2.5,0.8){\textcolor{red}{$\varepsilon_1$}}\psline[linecolor=red]{->}(2.5,0.8)(2,0.9)
%\uput[180](0.5,2.5){\textcolor{red}{$\varepsilon_{12}$}}\psline[linecolor=red]{->}(0.5,2.4)(1.7,1.6)\psline[linecolor=red]{->}(0.5,2.6)(1.4,4.3)
%\uput[0](1.8,3){\textcolor{red}{$\varepsilon_{23}$}}\psline[linecolor=red]{->}(2.3,2.7)(3.1,2.1)\psline[linecolor=red]{->}(2.3,3.3)(3.8,5)
%\uput[270](11.5,1){\textcolor{red}{$\tilde{\varepsilon}_1$}}\psline[linecolor=red]{->}(11.1,0.7)(8.9,1.8)\psline[linecolor=red]{->}(12,0.7)(13.7,1.2)
%\uput[180](7,5.8){\textcolor{red}{$\tilde{\varepsilon}_{23}$}}\psline[linecolor=red]{->}(7,5.6)(8.4,5)
\end{pspicture}
\caption{Branch loci of the three covers}\label{fig4}
\end{figure}

\begin{Remark}\rm
    Figures 3 and 4 are given for the case where $m_1<m_2<m_3$. If $m_1=m_2$ (resp. $m_2=m_3$) we have $\varepsilon_{12}=0$ (resp. $\varepsilon_{23}=0$) and it's easy to redraw the trees (which are then a little simpler). 
    In the very special case where $m_1=m_2=m_3$, the geometry of the branch locus is given in Figure 5.
    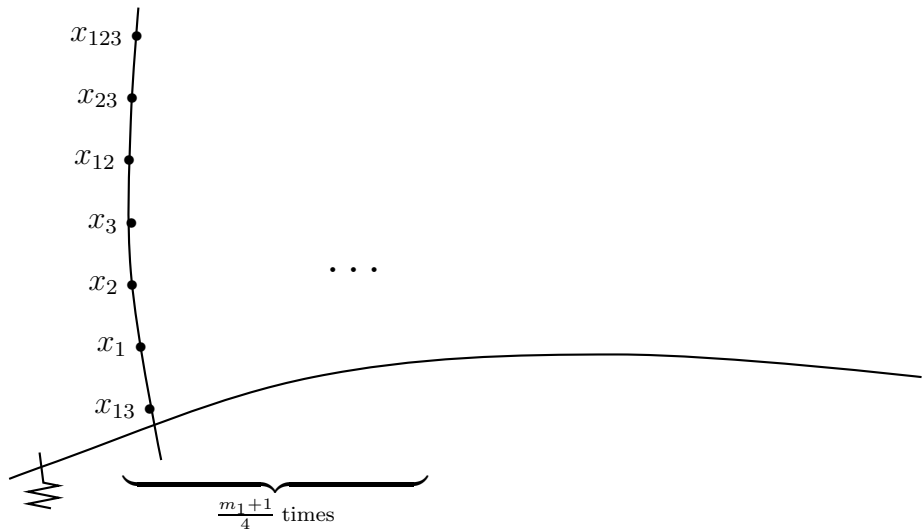
\begin{figure}[h]
\begin{pspicture}(0,-0.5)(16,7)	
%\pscurve[showpoints=true]{-}(0,1.3)
%(0.7,1.8)(3.3,0.5)(4,1.6)(0.4,0.4)
\psline[showpoints=false]{-}(0.4,0.6)(0.45,0.2)(0.65,0.16)(0.25,0.08)(0.65,0.0)(0.25,-0.08)(0.55,-0.14)%generique
\pscurve[showpoints=false]{-}(0,0.25)(0.4,0.4)(4,1.6)(8,1.9)(12,1.6)%horizontale0
\pscurve[showpoints=false]{-}(2,0.5)(1.6,3)(1.6,5)(1.7,6.5)%verticale1
%\pscurve[showpoints=false]{-}(0.5,4)(2.5,4.8)(4,5.2)(4.5,5.3)%horizontale1
%\pscurve[showpoints=false]{-}(9,1)(8.6,3.5)(8.6,5.5)(8.7,6.5)%verticale2

\psdot(1.73,2)\uput[180](1.73,2){$x_1$}
\psdot(1.58,4.48)\uput[180](1.58,4.48){$x_{12}$}
\psdot(1.62,2.82)\uput[180](1.62,2.82){$x_2$}
\psdot(1.61,3.64)\uput[180](1.61,3.64){$x_3$}
\psdot(1.85,1.18)\uput[180](1.85,1.18){$x_{13}$}
\psdot(1.62,5.3)\uput[180](1.62,5.3){$x_{23}$}
\psdot(1.68,6.12)\uput[180](1.68,6.12){$x_{123}$}

%\psdot(4,5.2)\uput[90](4,5.2){$x_{1}$}
%\psdot(8.6,3.5)\uput[0](8.6,3.5){$\hat{x}_{2}$}
%\psdot(8.6,5.5)\uput[0](8.6,5.5){$\check{x}_{2}$}

\uput[270](3.5,0.5){$\underbrace{\hspace{4cm}}_{\frac{m_1+1}{4}\;\mathrm{times}}$}
%\uput[270](10.5,0.5){$\underbrace{\hspace{4cm}}_{\frac{m_2-m_1}{2}\;\mathrm{times}}$}

\uput[0](4,3){\Large $\cdots$}%\uput[0](10,4){\Large $\cdots$}
%\psaxes{->}(0,0)(-0.5,-0.5)(14.5,7.5)[$x$,0][$y$, 90]	%creates axes
%\psdot(2,1)	%plots the point (2,1)
%\uput[0](2,1){$A$}	%labels the point (2,1) as A

\end{pspicture}
\caption{Branch locus when $m_1=m_2=m_3$}\label{fig5}
\end{figure}
    
   In the even more special case where $m_1+1=m_2+1=m_3+1=4$, we find equidistant geometry (which is consistent with Proposition 5.6 in ~\cite{Yang}). 
\end{Remark}
   
\subsection{Construction of covers with fixed geometry and reduction}
\subsubsection{Aim}

 The aim is as follows. We give $(r,s)\in\N^*\times\N$ ($\N$ is the set of non negative integers) and $X_1,\cdots,X_{r+s}\in R$, non-zero, such that $v(X_i-X_j)=0$ for all $i\neq j$. Write $n=2r+s$.  Let a $(\Z/2\Z)$-extension of $k[[t]]$ be given by an equation 
  $$\begin{aligned}
   w^2-w&=\sum_{\ell=0}^{n-1}\frac{a_{\ell}^2}{t^{2n-1-2\ell}}\end{aligned}\;\;\;(E')$$
  where $a_0,\cdots,a_{n-1}\in k$, $a_0\neq 0$. The conductor is then $2n$. 

\begin{figure}[h]
\begin{pspicture}(0,-0.5)(18,8)	
%\pscurve[showpoints=true]{-}(0,1.3)
%(0.7,1.8)(3.3,0.5)(4,1.6)(0.4,0.4)
\psline[showpoints=false]{-}(0.4,0.6)(0.45,0.2)(0.65,0.16)(0.25,0.08)(0.65,0.0)(0.25,-0.08)(0.55,-0.14)%generique
\pscurve[showpoints=false]{-}(0,0.25)(0.4,0.4)(4,1.6)(8,1.9)(12,1.6)(16,1)%horizontale0
\pscurve[showpoints=false]{-}(2,0.5)(1.6,3)(1.6,5)(1.7,6)%verticale1
%\pscurve[showpoints=false]{-}(2.4,2)(2.4,2.5)(2.5,3.5)(2.8,4.5)%verticale1bis
%\pscurve[showpoints=false]{-}(4,2.3)(4,2.8)(4.1,3.8)(4.4,4.8)%verticale1ter
%\pscurve[showpoints=false]{-}(3.2,1.6)(3.2,2.3)(3.3,3.3)(3.6,4.3)%verticale4
\pscurve[showpoints=false]{-}(0.5,1.5)(2.5,2.3)(4,2.7)(4.5,2.8)%horizontale1
\pscurve[showpoints=false]{-}(0.5,4)(2.5,4.8)(4,5.2)(4.5,5.3)%horizontale2
%\pscurve[showpoints=false]{-}(0.5,1)(2.5,1.8)(4,2.2)(4.5,2.3)%horizontale2
\pscurve[showpoints=false]{-}(9,1)(8.6,3.5)(8.6,5.5)(8.7,6.5)%verticale2

\uput[0](16,1){$I_0$}
\psdot(8.6,3.5)%\uput[0](8.6,3.5){$\tilde{x}_2$}
\psdot(8.6,5.5)%\uput[0](8.6,5.5){$\tilde{x}_{23}$}
%\uput[90](11,4.5){\textcolor{red}{$|\varrho_1^2\varrho_3|$}}\psline[linecolor=red]{->}(10.3,4.8)(8.9,4.3)

%\pscurve[showpoints=false]{-}(7.5,4.5)(9.5,5.3)(11,5.7)(11.5,5.8)%horizontale1
%\psdot(8.73,2.5)\uput[0](8.73,2.5){$\tilde{x}_3$}
%\psdot(9.5,5.3)\uput[90](9.5,5.3){$\tilde{x}_{2}$}
%\psdot(11,5.7)\uput[90](11,5.7){$\tilde{x}_{23}$}

%\psdot(2.3,1.72)\uput[90](2.3,1.72){$x_3$}
%\psdot(3.21,2.5)\uput[0](3.21,2.5){$x_{23}$}
%\psdot(3.33,3.5)\uput[0](3.33,3.5){$x_{2}$}
\psdot(2.5,2.3)%\uput[135](2.5,2.3){$x_{12}$}
\psdot(4,2.7)%\uput[135](4,2.7){$x_{123}$}
\psdot(2.5,4.8)%\uput[135](2.5,4.8){$x_{2}$}
\psdot(4,5.2)%\uput[135](4,5.2){$x_{23}$}

%\uput[0](5,4){\textcolor{red}{$|\varrho_3|$}}\psline[linecolor=red]{->}(5.3,4.3)(4.8,5.1)\psline[linecolor=red]{->}(5.3,3.6)(4.8,2.8)
%\uput[180](0,3){\textcolor{red}{$|\varrho_1|$}}\psline[linecolor=red]{->}(0,3)(1.2,3)
%\uput[90](16,3){\textcolor{red}{$1$}}\psline[linecolor=red]{->}(16,3)(15.5,1.3)

%\pscurve[showpoints=false]{-}(14,0.5)(13.6,3)(13.6,5)(13.7,6)%verticale1
%\psdot(13.6,3)\uput[0](13.6,3){$\check{x}_{3}$}
%\psdot(13.6,5)\uput[0](13.6,5){$\hat{x}_{3}$}
%\psdot(8.6,3.5)\uput[0](8.6,3.5){$\hat{x}_{2}$}
%\psdot(8.6,5.5)\uput[0](8.6,5.5){$\check{x}_{2}$}

\uput[270](3.5,0.5){$\underbrace{\hspace{4cm}}_{r\;\mathrm{times}}$}
\uput[270](10.5,0.5){$\underbrace{\hspace{3.5cm}}_{s\;\mathrm{times}}$}
%\uput[270](15,0.5){$\underbrace{\hspace{2.5cm}}_{\frac{m_3-m_2}{2}\;\mathrm{times}}$}

%\uput[0](7,3.5){\Large $\cdots$}%\uput[0](10,4){\Large $\cdots$}\uput[0](15,4){\Large $\cdots$}
%\psaxes{->}(0,0)(-0.5,-0.5)(14.5,7.5)[$x$,0][$y$, 90]	%creates axes
%\psdot(2,1)	%plots the point (2,1)
%\uput[0](2,1){$A$}	%labels the point (2,1) as A

%\uput[350](-1,1){\textcolor{red}{$\varepsilon_0$}}\psline[linecolor=red]{->}(-0.3,0.8)(0.3,0.5)
%\uput[0](2.5,0.8){\textcolor{red}{$\varepsilon_1$}}\psline[linecolor=red]{->}(2.5,0.8)(2,0.9)
%\uput[180](0.5,2.5){\textcolor{red}{$\varepsilon_{12}$}}\psline[linecolor=red]{->}(0.5,2.4)(1.7,1.6)\psline[linecolor=red]{->}(0.5,2.6)(1.4,4.3)
%\uput[0](1.8,3){\textcolor{red}{$\varepsilon_{23}$}}\psline[linecolor=red]{->}(2.3,2.7)(3.1,2.1)\psline[linecolor=red]{->}(2.3,3.3)(3.8,5)
%\uput[270](11.5,1){\textcolor{red}{$\tilde{\varepsilon}_1$}}\psline[linecolor=red]{->}(11.1,0.7)(8.9,1.8)\psline[linecolor=red]{->}(12,0.7)(13.7,1.2)
%\uput[180](7,5.8){\textcolor{red}{$\tilde{\varepsilon}_{23}$}}\psline[linecolor=red]{->}(7,5.6)(8.4,5)
\end{pspicture}
\caption{}\label{fig6}
$ $
\end{figure}   
 Finally, we give the following tree (see figure 6). We then want to find a polynomial $f\in R[X]$ such that 
  \begin{enumerate}[a.]
      \item The cover of the $2$-adic open disk given by the equation $Y^2=f(X)$ has good reduction relatively to a certain Gauss valuation, this reduction giving the equation $(E')$. %Pour vérifier cette bonne réduction, on utilisera le théorème 3.4 de ~\cite{GrMa}  : la bonne réduction est obtenue quand on a l'égalité des degrés des différentes générique et spéciale. Pour cela, il suffit de vérifier que le nombre de point de branchement à la fibre générique est égal au conducteur du revêtement d'Artin Schreier à la fibre spéciale.  
    \item The geometry of the branch locus is shown in Figure 6.
    \item The points $X_1,\cdots,X_{r+s}$ are part of the branch locus, i.e. $X_1,\cdots,X_{r+s}$ are simple roots of $f$.
  \end{enumerate}

As in the case  $(\Z/2\Z)^2$, the idea is to work from the exact differential form on the component $I_0$ of the tree, which is of the form 
  $$\frac{dx}{\prod\limits_{i=1}^r(x-x_i)^4\prod\limits_{j=1}^s(x-x_j)^2}$$
  and use its partial fraction expansion.

\subsubsection{Construction of the cover}

Let's begin with a lemma. 

\begin{lemma}
  Let $(r,s)\in\N^*\times\N$ and $X_1,\cdots,X_{r+s}\in R$, non-zero, such that $v(X_i-X_j)=0$ for all $i\neq j$. Let $n=2r+s$. Let $Q_1=\prod_{i=1}^r(X-X_i)$ and $Q_2=\prod_{j=r+1}^{r+s}(X-X_j)$ if $s>0$ and $Q_2=1$ if $s=0$. Let $\ell\in\cg0,n-1\cd$. Then there exist $\alpha_{i\ell},\beta_{i\ell},\gamma_{j\ell}\in R $ with $\alpha_{i\ell}$ and $\gamma_{j\ell}$ non-zero and such that  
  $$X^{2\ell}=Q_1^4(X)Q_2^2(X)\left(\sum_{i=1}^{r}\left(\frac{\alpha_{i\ell}^2}{(X-X_i)^4}+\frac{\beta_{i\ell}^2}{(X-X_i)^2}\right)+\sum_{j=r+1}^{r+s}\frac{\gamma_{j\ell}^2}{(X-X_j)^2}\right)+O(2) $$
  where $O(2)$ denotes a polynomial of $2R[X]$ of degree at most $2n-1$.
\end{lemma}

\begin{Proof}\rm
    We write the partial fraction expansion of $\frac{X^{\ell}}{Q_1^2Q_2}$. Then there exist $\alpha_{i\ell},\beta_{i\ell},\gamma_{j\ell}\in K:=Frac(R)$ such that 

     $$\frac{X^{\ell}}{Q_1^2(X)Q_2(X)}=\left(\sum_{i=1}^{r}\left(\frac{\alpha_{i\ell}}{(X-X_i)^2}+\frac{\beta_{i\ell}}{(X-X_i)}\right)+\sum_{j=r+1}^{r+s}\frac{\gamma_{j\ell}}{(X-X_j)}\right). $$

     We need to check that the coefficients $\alpha_{i\ell},\beta_{i\ell},\gamma_{j\ell}$ are indeed in $R$, but this is immediate since $v(X_i-X_j)=0$.

     We square this equality and multiply by $Q_1^4Q_2^2$, which completes the proof. 

     \begin{flushright}
     $\square$
     \end{flushright}
    
\end{Proof}

Let's fix some additional notations. Let $\rho_1,\rho_2\in R$ be such that $\frac{1}{2}v(2)\leq v(\rho_1)\leq v(\rho_2)$. Let $\rho=\rho_1^2\rho_2$ and assume that $v(\rho)<2v(2)$. Let $\rho_0=\left(\frac{4}{\rho}\right)^{\frac{1}{2n-1}}$. Note that $v(\rho_0)>0$.

Let $A_0, \cdots,A_{n-1}\in R$ such that $v(A_0)=0$. For all $i\in\cg1,r\cd$ and for $j\in\cg r+1,\cdots,r+s\cd$, write :

$$\begin{aligned}
  R_i(X)&=\sum_{\ell=0}^{n-1}\rho_0^{2\ell}A_{\ell}^2\left(\beta_{i\ell}^2(X-X_i)^3+\alpha_{i\ell}^2(X-X_i)\right)\\
  \tilde{R}_i(X)&=\sum_{\ell=0}^{n-1}\rho_0^{\ell}A_{\ell}\left(\beta_{i\ell}(X-X_i)^3+\alpha_{i\ell}(X-X_i)^2\right)X_i^{\frac{1}{2}}\\
  R_j(X)&=\sum_{\ell=0}^{n-1}\rho_0^{2\ell}A_{\ell}^2\gamma_{j\ell}^2(X-X_j)\\
   \tilde{R}_j(X)&=\sum_{\ell=0}^{n-1}\rho_0^{\ell}A_{\ell}\gamma_{j\ell}(X-X_j)X_j^{\frac{1}{2}}
\end{aligned}$$

Let $\delta_1,\cdots,\delta_r\in R$ and 
$$
\begin{aligned}
    \hat{Q}_i(X)&=(X-X_i)^2+\rho_1\delta_i(X-X_i)\\
    \hat{Q}(X)&=\prod_{i=1}^r \hat{Q}_i(X)\prod_{j=r+1}^{r+s}(X-X_j)\\
    P_i(X)&=\hat{Q}^2_i(X)+\rho R_i(X)+2\rho^{\frac{1}{2}}\tilde{R}_i(X)\\
     P_j(X)&=(X-X_j)^2+\rho R_j(X)+2\rho^{\frac{1}{2}}\tilde{R}_j(X)
\end{aligned}$$

Finally let $$f(X)=\prod_{i=1}^r P_i(X)\prod_{j=r+1}^{r+s}P_j(X).$$

\begin{Remark}\rm
    The $\delta_i$ may appear mysterious. In fact, they won't be involved in the calculations in section V.3. The idea is that we can choose them "as we please" without any influence on the equation we obtain in reduction. In the next part, we'll have to ensure that certain polynomials have roots in common, and we'll then have to adjust these coefficients. 
\end{Remark}
We then have the following lemma :
\begin{lemma}
    With the previous notations, we have the relation :
    \begin{align*}
\begin{split} 
f(X) &= \hat{Q}^2(X)+\rho Q_1^4(X)Q_2^2(X)\left(
\sum_{i=1}^{r}\frac{R_i(X)}{(X-X_i)^4}+\sum_{j=r+1}^{r+s}\frac{R_j(X)}{(X-X_j)^2}\right)
\\&+2\rho^{\frac{1}{2}}Q_1^4(X)Q_2^2(X)\left(
\sum_{i=1}^{r}\frac{\tilde{R}_i(X)}{(X-X_i)^4}+\sum_{j=r+1}^{r+s}\frac{\tilde{R}_j(X)}{(X-X_j)^2}\right)
\\&+O(4)
\end{split}
\end{align*}

where $O(4)$ denotes a polynomial in $4R[X]$ of degree at most $2n-1$.
\end{lemma}

\begin{Proof} \rm We have
\begin{align*}  
        f(X)&=\hat{Q}^2(X)\prod_{i=1}^r \frac{P_i(X)}{\hat{Q}_i^2(X)}\prod_{j=r+1}^{r+s}\frac{P_j(X)}{(X-X_j)^2}\\
        &=\hat{Q}^2(X)\prod_{i=1}^r \left[ 1+\rho\frac{R_i(X)}{\hat{Q}_i^2(X)}+2\rho^{\frac{1}{2}}\frac{\tilde{R}_i(X)}{\hat{Q}_i^2(X)}\right]
        \times \prod_{j=r+1}^{r+s} \left[ 1+\rho\frac{R_j(X)}{(X-X_j)^2}+2\rho^{\frac{1}{2}}\frac{\tilde{R}_j(X)}{(X-X_j)^2}\right]\\
        &=\hat{Q}^2(X)+\rho\hat{Q}^2(X)\left(\sum_{i=1}^r\frac{R_i(X)}{\hat{Q}_i^2(X)}+\sum_{j=r+1}^{r+s}\frac{R_j(X)}{(X-X_j)^2}\right)\\
        &+2\rho^{\frac{1}{2}}\hat{Q}^2(X)\left(\sum_{i=1}^r\frac{\tilde{R}_i(X)}{\hat{Q}_i^2(X)}+\sum_{j=r+1}^{r+s}\frac{\tilde{R}_j(X)}{(X-X_j)^2}\right)
        +O(4)
\end{align*}
since $v(\rho^2)\geq 2v(2)$ et $v(2\rho)\geq 2v(2)$.
 Or \begin{align*}
     \frac{\hat{Q}^2(X)}{\hat{Q}_i^2(X)}&=\prod_{\ell\neq i}\hat{Q}_{\ell}^2(X)\prod_{j=r+1}^{r+s}(X-X_j)^2\\
     &=\prod_{\ell\neq i} \left[(X-X_{\ell})^2+\rho_1\delta_{\ell}(X-X_{\ell})\right]^2\prod_{j=r+1}^{r+s}(X-X_j)^2\\
     &=\prod_{\ell\neq i}(X-X_{\ell})^4\prod_{j=r+1}^{r+s}(X-X_j)^2+O(\rho_1)\\
     &=\frac{Q_1^4(X)Q_2^2(X)}{(X-X_i)^4}+O(\rho_1)
 \end{align*} 

  In the same way, we show that 
  $$\frac{\hat{Q}^2(X)}{(X-X_j)^2}= \frac{Q_1^4(X)Q_2^2(X)}{(X-X_j)^2}+O(\rho_1).$$

  Given that $v(\rho\rho_1)\geq 2v(2)$, the equality of Lemma 4 is demonstrated.

   \begin{flushright}
     $\square$
     \end{flushright}
\end{Proof}

We can now prove the desired result.

\begin{Prop}
Let's take the same polynomial $f$ as for Lemma 4, and denote $a_0,\cdots,a_{n-1}$ the reductions modulo $\pi$ of $A_0,\cdots,A_{n-1}$. In particular, $a_0\neq 0$. Then the cover of $\mathbb{P}_K^1$ given by the equation $Y^2=f(X)$ has good reduction for the Gauss valuation relative to $T=\rho_0X$, the reduction being :
    $$w^2-w=\sum_{\ell=0}^{n-1}\frac{a_{\ell}^2}{t^{2n-1-2\ell}}.$$

   Furthermore, if $v(\rho_1)<v(\rho_2)$ and $v(\delta_i)=0$, the geometry of the branch locus is given by figure 7 (first tree). In the case where $v(\rho_1)=v(\rho_2)$, the geometry of the branch locus is given by figure 7 (second tree). 
\end{Prop}

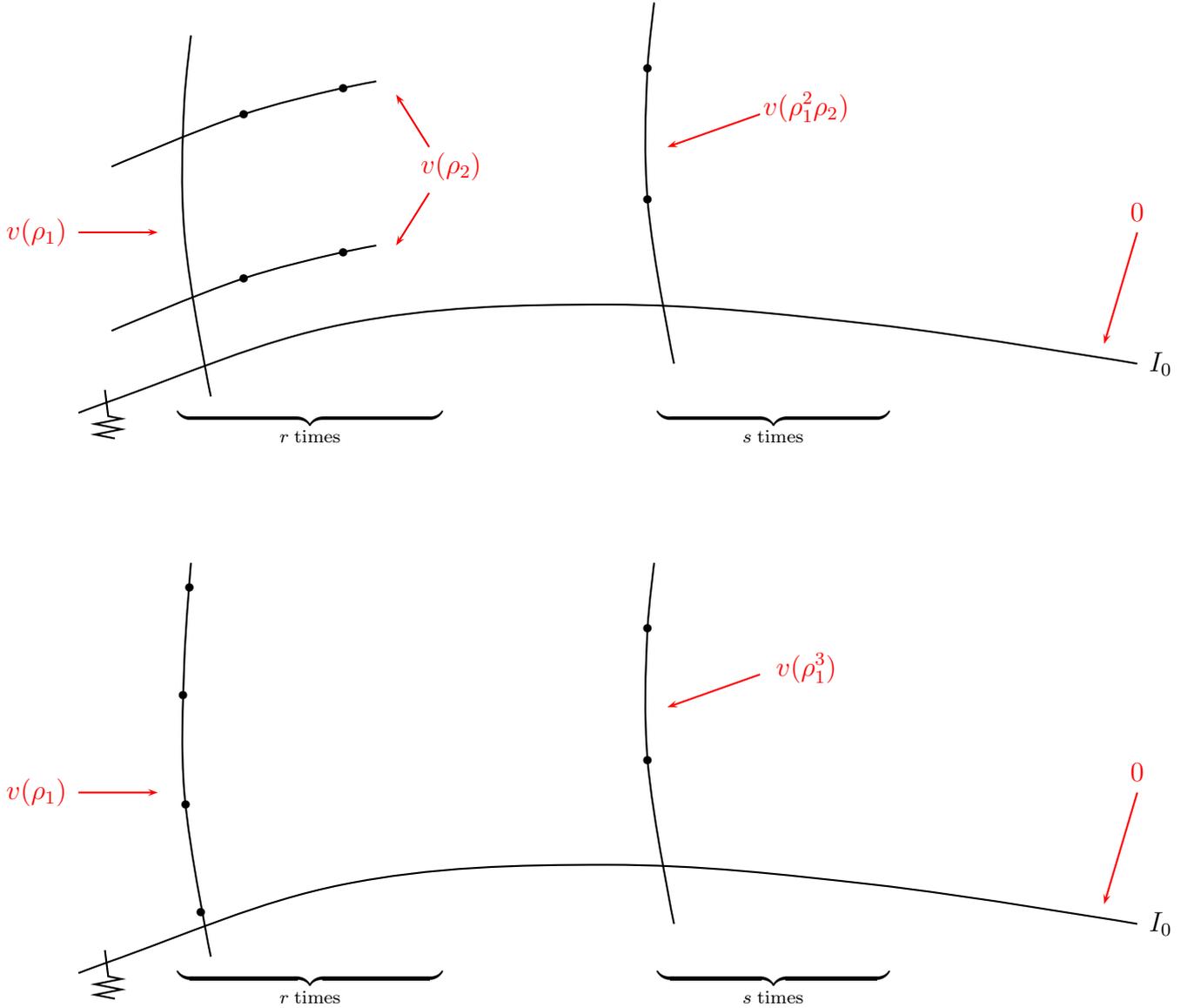
\begin{figure}[h]
\begin{pspicture}(0,-0.5)(18,6)	
%\pscurve[showpoints=true]{-}(0,1.3)
%(0.7,1.8)(3.3,0.5)(4,1.6)(0.4,0.4)
\psline[showpoints=false]{-}(0.4,0.6)(0.45,0.2)(0.65,0.16)(0.25,0.08)(0.65,0.0)(0.25,-0.08)(0.55,-0.14)%generique
\pscurve[showpoints=false]{-}(0,0.25)(0.4,0.4)(4,1.6)(8,1.9)(12,1.6)(16,1)%horizontale0
\pscurve[showpoints=false]{-}(2,0.5)(1.6,3)(1.6,5)(1.7,6)%verticale1
%\pscurve[showpoints=false]{-}(2.4,2)(2.4,2.5)(2.5,3.5)(2.8,4.5)%verticale1bis
%\pscurve[showpoints=false]{-}(4,2.3)(4,2.8)(4.1,3.8)(4.4,4.8)%verticale1ter
%\pscurve[showpoints=false]{-}(3.2,1.6)(3.2,2.3)(3.3,3.3)(3.6,4.3)%verticale4
\pscurve[showpoints=false]{-}(0.5,1.5)(2.5,2.3)(4,2.7)(4.5,2.8)%horizontale1
\pscurve[showpoints=false]{-}(0.5,4)(2.5,4.8)(4,5.2)(4.5,5.3)%horizontale2
%\pscurve[showpoints=false]{-}(0.5,1)(2.5,1.8)(4,2.2)(4.5,2.3)%horizontale2
\pscurve[showpoints=false]{-}(9,1)(8.6,3.5)(8.6,5.5)(8.7,6.5)%verticale2

\uput[0](16,1){$I_0$}
\psdot(8.6,3.5)%\uput[0](8.6,3.5){$\tilde{x}_2$}
\psdot(8.6,5.5)%\uput[0](8.6,5.5){$\tilde{x}_{23}$}
\uput[90](11,4.5){\textcolor{red}{$v(\rho_1^2\rho_2)$}}\psline[linecolor=red]{->}(10.3,4.8)(8.9,4.3)

%\pscurve[showpoints=false]{-}(7.5,4.5)(9.5,5.3)(11,5.7)(11.5,5.8)%horizontale1
%\psdot(8.73,2.5)\uput[0](8.73,2.5){$\tilde{x}_3$}
%\psdot(9.5,5.3)\uput[90](9.5,5.3){$\tilde{x}_{2}$}
%\psdot(11,5.7)\uput[90](11,5.7){$\tilde{x}_{23}$}

%\psdot(2.3,1.72)\uput[90](2.3,1.72){$x_3$}
%\psdot(3.21,2.5)\uput[0](3.21,2.5){$x_{23}$}
%\psdot(3.33,3.5)\uput[0](3.33,3.5){$x_{2}$}
\psdot(2.5,2.3)%\uput[135](2.5,2.3){$x_{12}$}
\psdot(4,2.7)%\uput[135](4,2.7){$x_{123}$}
\psdot(2.5,4.8)%\uput[135](2.5,4.8){$x_{2}$}
\psdot(4,5.2)%\uput[135](4,5.2){$x_{23}$}

\uput[0](5,4){\textcolor{red}{$v(\rho_2)$}}\psline[linecolor=red]{->}(5.3,4.3)(4.8,5.1)\psline[linecolor=red]{->}(5.3,3.6)(4.8,2.8)
\uput[180](0,3){\textcolor{red}{$v(\rho_1)$}}\psline[linecolor=red]{->}(0,3)(1.2,3)
\uput[90](16,3){\textcolor{red}{$0$}}\psline[linecolor=red]{->}(16,3)(15.5,1.3)

%\pscurve[showpoints=false]{-}(14,0.5)(13.6,3)(13.6,5)(13.7,6)%verticale1
%\psdot(13.6,3)\uput[0](13.6,3){$\check{x}_{3}$}
%\psdot(13.6,5)\uput[0](13.6,5){$\hat{x}_{3}$}
%\psdot(8.6,3.5)\uput[0](8.6,3.5){$\hat{x}_{2}$}
%\psdot(8.6,5.5)\uput[0](8.6,5.5){$\check{x}_{2}$}

\uput[270](3.5,0.5){$\underbrace{\hspace{4cm}}_{r\;\mathrm{times}}$}
\uput[270](10.5,0.5){$\underbrace{\hspace{3.5cm}}_{s\;\mathrm{times}}$}
%\uput[270](15,0.5){$\underbrace{\hspace{2.5cm}}_{\frac{m_3-m_2}{2}\;\mathrm{times}}$}

%\uput[0](7,3.5){\Large $\cdots$}%\uput[0](10,4){\Large $\cdots$}\uput[0](15,4){\Large $\cdots$}
%\psaxes{->}(0,0)(-0.5,-0.5)(14.5,7.5)[$x$,0][$y$, 90]	%creates axes
%\psdot(2,1)	%plots the point (2,1)
%\uput[0](2,1){$A$}	%labels the point (2,1) as A

%\uput[350](-1,1){\textcolor{red}{$\varepsilon_0$}}\psline[linecolor=red]{->}(-0.3,0.8)(0.3,0.5)
%\uput[0](2.5,0.8){\textcolor{red}{$\varepsilon_1$}}\psline[linecolor=red]{->}(2.5,0.8)(2,0.9)
%\uput[180](0.5,2.5){\textcolor{red}{$\varepsilon_{12}$}}\psline[linecolor=red]{->}(0.5,2.4)(1.7,1.6)\psline[linecolor=red]{->}(0.5,2.6)(1.4,4.3)
%\uput[0](1.8,3){\textcolor{red}{$\varepsilon_{23}$}}\psline[linecolor=red]{->}(2.3,2.7)(3.1,2.1)\psline[linecolor=red]{->}(2.3,3.3)(3.8,5)
%\uput[270](11.5,1){\textcolor{red}{$\tilde{\varepsilon}_1$}}\psline[linecolor=red]{->}(11.1,0.7)(8.9,1.8)\psline[linecolor=red]{->}(12,0.7)(13.7,1.2)
%\uput[180](7,5.8){\textcolor{red}{$\tilde{\varepsilon}_{23}$}}\psline[linecolor=red]{->}(7,5.6)(8.4,5)
\end{pspicture}

\begin{pspicture}(0,-0.5)(18,8)	
%\pscurve[showpoints=true]{-}(0,1.3)
%(0.7,1.8)(3.3,0.5)(4,1.6)(0.4,0.4)
\psline[showpoints=false]{-}(0.4,0.6)(0.45,0.2)(0.65,0.16)(0.25,0.08)(0.65,0.0)(0.25,-0.08)(0.55,-0.14)%generique
\pscurve[showpoints=false]{-}(0,0.25)(0.4,0.4)(4,1.6)(8,1.9)(12,1.6)(16,1)%horizontale0
\pscurve[showpoints=false]{-}(2,0.5)(1.6,3)(1.6,5)(1.7,6.5)%verticale1
%\pscurve[showpoints=false]{-}(2.4,2)(2.4,2.5)(2.5,3.5)(2.8,4.5)%verticale1bis
%\pscurve[showpoints=false]{-}(4,2.3)(4,2.8)(4.1,3.8)(4.4,4.8)%verticale1ter
%\pscurve[showpoints=false]{-}(3.2,1.6)(3.2,2.3)(3.3,3.3)(3.6,4.3)%verticale4
%\pscurve[showpoints=false]{-}(0.5,1.5)(2.5,2.3)(4,2.7)(4.5,2.8)%horizontale1
%\pscurve[showpoints=false]{-}(0.5,4)(2.5,4.8)(4,5.2)(4.5,5.3)%horizontale2
%\pscurve[showpoints=false]{-}(0.5,1)(2.5,1.8)(4,2.2)(4.5,2.3)%horizontale2
\pscurve[showpoints=false]{-}(9,1)(8.6,3.5)(8.6,5.5)(8.7,6.5)%verticale2

\uput[0](16,1){$I_0$}
\psdot(8.6,3.5)%\uput[0](8.6,3.5){$\tilde{x}_2$}
\psdot(8.6,5.5)%\uput[0](8.6,5.5){$\tilde{x}_{23}$}
\uput[90](11,4.5){\textcolor{red}{$v(\rho_1^3)$}}\psline[linecolor=red]{->}(10.3,4.8)(8.9,4.3)

%\psdot(1.73,2)\uput[180](1.73,2){$x_1$}
\psdot(1.58,4.48)%\uput[180](1.58,4.48){$x_{12}$}
\psdot(1.62,2.82)%\uput[180](1.62,2.82){$x_2$}
%\psdot(1.61,3.64)\uput[180](1.61,3.64){$x_3$}
\psdot(1.85,1.18)%\uput[180](1.85,1.18){$x_{13}$}
%\psdot(1.62,5.3)\uput[180](1.62,5.3){$x_{23}$}
\psdot(1.68,6.12)%\uput[180](1.68,6.12){$x_{123}$}

%\pscurve[showpoints=false]{-}(7.5,4.5)(9.5,5.3)(11,5.7)(11.5,5.8)%horizontale1
%\psdot(8.73,2.5)\uput[0](8.73,2.5){$\tilde{x}_3$}
%\psdot(9.5,5.3)\uput[90](9.5,5.3){$\tilde{x}_{2}$}
%\psdot(11,5.7)\uput[90](11,5.7){$\tilde{x}_{23}$}

%\psdot(2.3,1.72)\uput[90](2.3,1.72){$x_3$}
%\psdot(3.21,2.5)\uput[0](3.21,2.5){$x_{23}$}
%\psdot(3.33,3.5)\uput[0](3.33,3.5){$x_{2}$}
%\psdot(2.5,2.3)%\uput[135](2.5,2.3){$x_{12}$}
%\psdot(4,2.7)%\uput[135](4,2.7){$x_{123}$}
%\psdot(2.5,4.8)%\uput[135](2.5,4.8){$x_{2}$}
%\psdot(4,5.2)%\uput[135](4,5.2){$x_{23}$}

%\uput[0](5,4){\textcolor{red}{$|\rho_2|$}}\psline[linecolor=red]{->}(5.3,4.3)(4.8,5.1)\psline[linecolor=red]{->}(5.3,3.6)(4.8,2.8)
\uput[180](0,3){\textcolor{red}{$v(\rho_1)$}}\psline[linecolor=red]{->}(0,3)(1.2,3)
\uput[90](16,3){\textcolor{red}{$0$}}\psline[linecolor=red]{->}(16,3)(15.5,1.3)

%\pscurve[showpoints=false]{-}(14,0.5)(13.6,3)(13.6,5)(13.7,6)%verticale1
%\psdot(13.6,3)\uput[0](13.6,3){$\check{x}_{3}$}
%\psdot(13.6,5)\uput[0](13.6,5){$\hat{x}_{3}$}
%\psdot(8.6,3.5)\uput[0](8.6,3.5){$\hat{x}_{2}$}
%\psdot(8.6,5.5)\uput[0](8.6,5.5){$\check{x}_{2}$}

\uput[270](3.5,0.5){$\underbrace{\hspace{4cm}}_{r\;\mathrm{times}}$}
\uput[270](10.5,0.5){$\underbrace{\hspace{3.5cm}}_{s\;\mathrm{times}}$}
%\uput[270](15,0.5){$\underbrace{\hspace{2.5cm}}_{\frac{m_3-m_2}{2}\;\mathrm{times}}$}

%\uput[0](7,3.5){\Large $\cdots$}%\uput[0](10,4){\Large $\cdots$}\uput[0](15,4){\Large $\cdots$}
%\psaxes{->}(0,0)(-0.5,-0.5)(14.5,7.5)[$x$,0][$y$, 90]	%creates axes
%\psdot(2,1)	%plots the point (2,1)
%\uput[0](2,1){$A$}	%labels the point (2,1) as A

%\uput[350](-1,1){\textcolor{red}{$\varepsilon_0$}}\psline[linecolor=red]{->}(-0.3,0.8)(0.3,0.5)
%\uput[0](2.5,0.8){\textcolor{red}{$\varepsilon_1$}}\psline[linecolor=red]{->}(2.5,0.8)(2,0.9)
%\uput[180](0.5,2.5){\textcolor{red}{$\varepsilon_{12}$}}\psline[linecolor=red]{->}(0.5,2.4)(1.7,1.6)\psline[linecolor=red]{->}(0.5,2.6)(1.4,4.3)
%\uput[0](1.8,3){\textcolor{red}{$\varepsilon_{23}$}}\psline[linecolor=red]{->}(2.3,2.7)(3.1,2.1)\psline[linecolor=red]{->}(2.3,3.3)(3.8,5)
%\uput[270](11.5,1){\textcolor{red}{$\tilde{\varepsilon}_1$}}\psline[linecolor=red]{->}(11.1,0.7)(8.9,1.8)\psline[linecolor=red]{->}(12,0.7)(13.7,1.2)
%\uput[180](7,5.8){\textcolor{red}{$\tilde{\varepsilon}_{23}$}}\psline[linecolor=red]{->}(7,5.6)(8.4,5)
\end{pspicture}
\caption{Geometry in the case $v(\rho_1)<v(\rho_2),v(\delta_i)=0$ and in the case $v(\rho_1)=v(\rho_2)$ }\label{fig7}
$ $
\end{figure}   

\begin{Remark}\rm
    The case $v(\rho_1)<v(\rho_2),v(\delta_i)=0$ is the one that will appear in the most general case where $m_1<m_2<m_3$. The case $v(\rho_1)=v(\rho_2)$ arises in the special case where two conductors are equal. For example, in the case where $m_1=m_2<m_3$, $Rev_1$ and $Rev_2$ will have a geometry corresponding to figure 7 (first tree) and $Rev_3$ to figure 7 (second tree). 
\end{Remark}

\begin{Proof}\rm
    Let's start with the result of Lemma 2 and look separately at the $\rho$ and $2\rho^{\frac{1}{2}}$ terms.  We have
    \begin{align*}
        &\rho Q_1^4(X)Q_2^2(X)\left(
\sum_{i=1}^{r}\frac{R_i(X)}{(X-X_i)^4}+\sum_{j=r+1}^{r+s}\frac{R_j(X)}{(X-X_j)^2}\right)\\
=&\rho Q_1^4(X)Q_2^2(X)\sum_{\ell=0}^{n-1}\rho_0^{2\ell}A_{\ell}^2
\left(\sum_{i=1}^{r}\frac{\alpha_{i\ell}^2}{(X-X_i)^3}+\frac{\beta_{i\ell}^2}{(X-X_i)}+\sum_{j=r+1}^{r+s}\frac{\gamma_{j\ell}^2}{(X-X_j)}\right)\\
=&\rho Q_1^4(X)Q_2^2(X)\sum_{\ell=0}^{n-1}\rho_0^{2\ell}A_{\ell}^2
\left(\sum_{i=1}^{r}\frac{\alpha_{i\ell}^2(X-X_i)}{(X-X_i)^4}+\frac{\beta_{i\ell}^2(X-X_i)}{(X-X_i)^2}+\sum_{j=r+1}^{r+s}\frac{\gamma_{j\ell}^2(X-X_j)}{(X-X_j)^2}\right)\\
&=\rho Q_1^4(X)Q_2^2(X)\sum_{\ell=0}^{n-1}\rho_0^{2\ell}A_{\ell}^2
\left(\sum_{i=1}^{r}\frac{\alpha_{i\ell}^2}{(X-X_i)^4}+\frac{\beta_{i\ell}^2}{(X-X_i)^2}+\sum_{j=r+1}^{r+s}\frac{\gamma_{j\ell}^2}{(X-X_j)^2}\right)X\\
&-\rho Q_1^4(X)Q_2^2(X)\sum_{\ell=0}^{n-1}\rho_0^{2\ell}A_{\ell}^2
\left(\sum_{i=1}^{r}\frac{\alpha_{i\ell}^2X_i}{(X-X_i)^4}+\frac{\beta_{i\ell}^2X_i}{(X-X_i)^2}+\sum_{j=r+1}^{r+s}\frac{\gamma_{j\ell}^2X_j}{(X-X_j)^2}\right)\\
&=\rho Q_1^4(X)Q_2^2(X)\sum_{\ell=0}^{n-1}\rho_0^{2\ell}A_{\ell}^2\left(
\sum_{i=1}^{r}\frac{\alpha_{i\ell}X_i^{\frac{1}{2}}}{(X-X_i)^2}+\frac{\beta_{i\ell}X_i^{\frac{1}{2}}}{(X-X_i)}+\sum_{j=r+1}^{r+s}\frac{\gamma_{j\ell}X_j^{\frac{1}{2}}}{(X-X_j)}\right)^2\\&+\rho\sum_{\ell=0}^{n-1}\rho_0^{2\ell}A_{\ell}^2X^{2\ell+1}+O(4)
    \end{align*}

   taking into account Lemma 3 and the fact that $v(2\rho)>2v(2)$. 

    On the other hand, we have 
    \begin{align*}
        &2\rho^{\frac{1}{2}}Q_1^4(X)Q_2^2(X)\left(
\sum_{i=1}^{r}\frac{\tilde{R}_i(X)}{(X-X_i)^4}+\sum_{j=r+1}^{r+s}\frac{\tilde{R}_j(X)}{(X-X_j)^2}\right)\\
=&2\rho^{\frac{1}{2}}Q_1^4(X)Q_2^2(X)\sum_{\ell=0}^{n-1}\rho_0^{\ell}A_{\ell}\left(
\sum_{i=1}^{r}\frac{\alpha_{i\ell}X_i^{\frac{1}{2}}}{(X-X_i)^2}+\frac{\beta_{i\ell}X_i^{\frac{1}{2}}}{(X-X_i)}+\sum_{j=r+1}^{r+s}\frac{\gamma_{j\ell}X_j^{\frac{1}{2}}}{(X-X_j)}\right)\\
=&2\rho^{\frac{1}{2}}Q_1^2(X)Q_2(X)\hat{Q}(X)\sum_{\ell=0}^{n-1}\rho_0^{\ell}A_{\ell}\left(
\sum_{i=1}^{r}\frac{\alpha_{i\ell}X_i^{\frac{1}{2}}}{(X-X_i)^2}+\frac{\beta_{i\ell}X_i^{\frac{1}{2}}}{(X-X_i)}+\sum_{j=r+1}^{r+s}\frac{\gamma_{j\ell}X_j^{\frac{1}{2}}}{(X-X_j)}\right)+O(4)
    \end{align*}

since $\hat{Q}(X)=Q_1^2(X)Q_2(X)+O(\rho_1)$ and $v(2\rho^{\frac{1}{2}}\rho_1)\geq2v(2)$. 

Thus (using the result of Lemma 4)
\begin{align*}
    f(X)&=\hat{Q}^2(X)+\rho\sum_{\ell=0}^{n-1}\rho_0^{2\ell}A_{\ell}^2X^{2\ell+1}\\&+\rho Q_1^4(X)Q_2^2(X)\sum_{\ell=0}^{n-1}\rho_0^{2\ell}A_{\ell}^2\left(
\sum_{i=1}^{r}\frac{\alpha_{i\ell}X_i^{\frac{1}{2}}}{(X-X_i)^2}+\frac{\beta_{i\ell}X_i^{\frac{1}{2}}}{(X-X_i)}+\sum_{j=r+1}^{r+s}\frac{\gamma_{j\ell}X_j^{\frac{1}{2}}}{(X-X_j)}\right)^2\\
&+2\rho^{\frac{1}{2}}Q_1^2(X)Q_2(X)\hat{Q}(X)\sum_{\ell=0}^{n-1}\rho_0^{\ell}A_{\ell}\left(
\sum_{i=1}^{r}\frac{\alpha_{i\ell}X_i^{\frac{1}{2}}}{(X-X_i)^2}+\frac{\beta_{i\ell}X_i^{\frac{1}{2}}}{(X-X_i)}+\sum_{j=r+1}^{r+s}\frac{\gamma_{j\ell}X_j^{\frac{1}{2}}}{(X-X_j)}\right)+O(4)\\
&=\left[\hat{Q}(X)+\rho^{\frac{1}{2}}Q_1^2(X)Q_2(X)\sum_{\ell=0}^{n-1}\rho_0^{\ell}A_{\ell}\left(
\sum_{i=1}^{r}\frac{\alpha_{i\ell}X_i^{\frac{1}{2}}}{(X-X_i)^2}+\frac{\beta_{i\ell}X_i^{\frac{1}{2}}}{(X-X_i)}+\sum_{j=r+1}^{r+s}\frac{\gamma_{j\ell}X_j^{\frac{1}{2}}}{(X-X_j)}\right)\right]^2\\&+\rho\sum_{\ell=0}^{n-1}\rho_0^{2\ell}A_{\ell}^2X^{2\ell+1}+O(4)\\
&=Q_0^2(X)+\rho\sum_{\ell=0}^{n-1}\rho_0^{2\ell}A_{\ell}^2X^{2\ell+1}+O(4)
\end{align*}
by setting $Q_0(X)=\hat{Q}(X)+\rho^{\frac{1}{2}}Q_1^2(X)Q_2(X)\sum_{\ell=0}^{n-1}\rho_0^{\ell}A_{\ell}\left(
\sum_{i=1}^{r}\frac{\alpha_{i\ell}X_i^{\frac{1}{2}}}{(X-X_i)^2}+\frac{\beta_{i\ell}X_i^{\frac{1}{2}}}{(X-X_i)}+\sum_{j=r+1}^{r+s}\frac{\gamma_{j\ell}X_j^{\frac{1}{2}}}{(X-X_j)}\right)$.

  Let $Y=2V+Q_0(X)$. The equation $Y^2=f(X)$ then becomes
  $$4(V^2+VQ_0(X))= \rho\sum_{\ell=0}^{n-1}\rho_0^{2\ell}A_{\ell}^2X^{2\ell+1}+O(4) $$ i.e.
  $$ \frac{4(V^2+VQ_0(X))}{X^{2n}}=\rho\rho_0^{2n-1}\sum_{\ell=0}^{n-1}A_{\ell}^2\left(\frac{1}{\rho_0X}\right)^{2n-1-2\ell}+\frac{O(4)}{X^{2n}}.$$

 Now $\rho\rho_0^{2n-1}=4$, so we get
  $$ \frac{(V^2+VQ_0(X))}{X^{2n}}=\sum_{\ell=0}^{n-1}A_{\ell}^2\left(\frac{1}{\rho_0X}\right)^{2n-1-2\ell}+\frac{O(1)}{X^{2n}}$$

  Here $O(1)$ denotes an element of $R[X]$ of degree at most $2n-1$. 
  Let's write $T=\rho_0X$ and $V=X^nW$. The equation becomes
  $$W^2+W\frac{Q_0(X)}{X^{n}}=\sum_{\ell=0}^{n-1}A_{\ell}^2\left(\frac{1}{T}\right)^{2n-1-2\ell}+\frac{O(1)}{X^{2n}}, $$
  which in reduction gives
  $$w^2-w=\sum_{\ell=0}^{n-1}\frac{a_{\ell}^2}{t^{2n-1-2\ell}}.$$

Let's take a look at the geometry of the branch locus. First of all, it's clear that for any $i\in\cg1,r+s\cd $, the roots of $P_i$ are simple (a little discriminant calculation shows this easily). Moreover, the roots of $P_i$ belong to the open disk of center $X_i$ and zero valuation. Since we have $v(X_i-X_j)=0$ as soon as $i\neq j$, it follows that a root of $P_i$ and a root of $P_j$ are distant by 1 as soon as $i\neq j$. 
%This justifies that the same applies to the roots of $P_i$ and $P_j$.

Let's fix a $i\in\cg1,r\cd$. To see that the geometry of the four roots of $P_i$ is as announced, simply note that $X_i$ is a root of $P_i$ and that, for any other root $Y_i$ of $P_i$, we can easily calculate $v(X_i-Y_i)$ using Newton's polygon of $\frac{P_i}{X-X_i}$. 

 (For readers unfamiliar with the concept of Newton polygon : the Newton polygon of a polynomial $P=\sum_{k=0}^na_kX^k$ with coefficients in $K$ is defined to be the lower boundary of the convex hull of the set of points $\left(i,v(a_i)\right)$ ignoring the points with $a_i = 0$. The knowledge of the slopes of the line segments of the Newton polygon then gives the valuation of the roots of $P$.)

   \begin{flushright}
     $\square$
     \end{flushright}
\end{Proof}

\begin{Remark}\rm
    We preserve the good reduction (and the same equation in reduction) if we replace $f(X)$ by $f(X)+O(4)$, where $O(4)\in4R[X]$ is a polynomial of degree at most $2n-1$. 
\end{Remark}

\subsection{A slightly simplified version of proposition 2}
We'll state a similar result by slightly modifying the polynomials $P_i$. This is the subject of proposition 3, and it's this proposition that will be useful to us in what follows.

We can write 
$$P_i(X)=\left[(X-X_i)^2+\rho_1\delta_i(X-X_i)\right]^2+\rho(\beta_i(X-X_i)^3+\alpha_i(X-X_i))+2\rho^{\frac{1}{2}}\mu_i(X-X_i)^2$$
setting
$$\alpha_i=\sum_{\ell=0}^{n-1}\rho_0^{2\ell}A_{\ell}^2\alpha_{i\ell}^2$$
$$\beta_i=\sum_{\ell=0}^{n-1}\rho_0^{2\ell}A_{\ell}^2\beta_{i\ell}^2+\frac{2}{\rho^{\frac{1}{2}}}\rho_0^{\ell}A_{\ell}\beta_{i\ell}X_i^{\frac{1}{2}}$$
$$\mu_i=\sum_{\ell=0}^{n-1}\rho_0^{\ell}A_{\ell}\alpha_{i\ell}X_i^{\frac{1}{2}}.$$
Similarly, we write $P_j(X)=(X-X_j)^2+\rho\gamma_j(X-T_j)$ where

$$\gamma_j=\sum_{\ell=0}^{n-1}\rho_0^{2\ell}A_{\ell}^2\gamma_{j\ell}^2+\frac{2}{\rho^{\frac{1}{2}}}\rho_0^{\ell}A_{\ell}\gamma_{j\ell}X_j^{\frac{1}{2}}$$

The idea here is to eliminate the last term in $(X-X_i)^2$. We write
\begin{align*}
    P_i&=\left[(X-X_i)^2+\rho_1\delta_i(X-X_i)+2^{\frac{1}{2}}\rho^{\frac{1}{4}}\mu_i^{\frac{1}{2}}(X-X_i)\right]^2
\\&+(\rho\beta_i-2^{\frac{3}{2}}\rho^{\frac{1}{4}}\mu_i^{\frac{1}{2}})(X-X_i)^3+\rho\alpha_i(X-X_i))+O(4)
\end{align*}
where $O(4)\in 4R[X]$ is of degree 3 or less.

Thus setting 
$$\beta'_i=\beta_i-2^{\frac{3}{2}}\rho^{\frac{-3}{4}}\mu_i^{\frac{1}{2}}$$
$$\delta'_i=\delta_i+2^{\frac{1}{2}}\frac{\rho^{\frac{1}{4}}}{\rho_1}\mu_i^{\frac{1}{2}}$$ and
 $$ \tilde{P}_i(X)=\left[(X-X_i)^2+\rho_1\delta'_i(X-X_i)\right]^2+\rho\beta'_i(X-X_i)^3+\rho\alpha_i(X-X_i) $$

 so that $\tilde{P}_i(X)=P_i(X)+O(4)$. We have then the
 \begin{Prop}
     Let $\tilde{f}(X)=\prod_{i=1}^r \tilde{P}_i(X)\prod_{j=r+1}^{r+s}P_j(X)$. Then the cover of $\mathbb{P}_K^1$ given by the equation $Y^2=\tilde{f}(X)$ has good reduction for the Gauss valuation relative to $T=\rho_0X$, the reduction being :
    $$w^2-w=\sum_{\ell=0}^{n-1}\frac{a_{\ell}^2}{t^{2n-1-2\ell}}.$$
 \end{Prop}

 \begin{Proof}\rm
    We have $\tilde{f}(X)=f(X)+O(4)$ where $O(4)\in 4R[X]$ is of degree at most $2n-1$. It is therefore sufficient to refer to Remark 4.

     \begin{flushright}
     $\square$
     \end{flushright}
 \end{Proof}

 \begin{Remark}\rm
   Since the $\delta_i$ could be chosen as desired in Proposition 2, it follows that the $\delta'_i$ can be chosen as desired in Proposition 3.
 \end{Remark}
%\section{Le théorème}
\subsection{Application of proposition 3}
Consider a $(Z/2\Z)^3$-extension of $k[[t]]$ given by the equations 
$$\left\{
\begin{aligned}
    w_1^2-w_1&=\sum_{\ell=0}^{\frac{m_1-1}{2}}\frac{a_{\ell}^2}{t^{m_1-2\ell}}\;\;\;\;(1)\\
    w_2^2-w_2&=\sum_{\ell=0}^{\frac{m_2-1}{2}}\frac{b_{\ell}^2}{t^{m_2-2\ell}}\;\;\;\;(2)\\
    w_3^2-w_3&=\sum_{\ell=0}^{\frac{m_3-1}{2}}\frac{c_{\ell}^2}{t^{m_3-2\ell}}\;\;\;\;(3)
\end{aligned}
\right.$$
where $a_{\ell}$, $b_{\ell}$, $c_{\ell} \in k$ et $a_0,b_0,c_0\in k^*$.

Recall that the triplet $(m_1,m_2,m_3)$ is minimal (in the sense of definition 1) and that $4$ divides $m_1+1$. So, in the case where $m_1=m_2$ (for example), the elements $a_0^2$ and $b_0^2$ are $\F_2$-independent (and therefore $a_0^2+b_0^2\neq 0$).

The idea is to use Proposition 3 to build three covers which have simultaneously good reduction (the reductions being those given by the equations above)  and such that the number of common branch points of the three covers is exactly $\frac{m_1+1}{4}$.

Recall also the notations :
$$m=2m_3+m_2-m_1,\;\;\varrho_1=2^{\frac{1}{2}},\;\;
\varrho_2=2^{\frac{m_3+m_2-m_1}{m}},\;\;
\varrho_3=2^{\frac{2m_3-m_1}{m}}.$$

and $v(\varrho_1)\leq v(\varrho_2)\leq v(\varrho_3)$.

Consider $X_1,\cdots, X_{\frac{m_1+1}{4}+\frac{m_3-m_1}{2}}\in R$, non-zero and such that $v(X_i-X_j)=0$ for $i\neq j$. We'll apply Proposition 3 three times, with judicious choices for $(r,s)$, $(\rho_1,\rho_2)$ and $X_1,\cdots, X_{r+s}$.
\begin{enumerate}[$\bullet$]
    \item Let's first apply Proposition 3 to the case where $r=\frac{m_1+1}{4}$, $s=0$ , $(\rho_1,\rho_2)=(\varrho_2,\varrho_3)$  and $A_{\ell}$ lifting $a_{\ell}$. To do this, we easily check that the conditions $\frac{1}{2}v(2)\leq v(\varrho_2)\leq v(\varrho_3)$ and $v(\varrho_2^2\varrho_3)<2v(2)$ are satisfied. This gives  $n=2r=\frac{m_1+1}{2}$ and $\rho=\varrho_2^2\varrho_3$. We then find a cover of equation $Y^2=f_1(X)$ where $f_1(X)=\prod_{i=1}^{\frac{m_1+1}{4}}\tilde{P}_{i,1}(X)$ and
    $$\tilde{P}_{i,1}(X)=\left[ (X-X_i)^2+\varrho_2\delta'_{i,1}(X-X_i)\right]^2+\varrho_2^2\varrho_3\beta'_{i,1}(X-X_i)^3+\varrho_2^2\varrho_3\alpha_{i,1}(X-X_i).$$
    In this case we have 
    $$ \rho_0=\left(\frac{4}{\rho}\right)^{\frac{1}{4r-1}}=
    \left(\frac{4}{\varrho_2^2\varrho_3}\right)^{\frac{1}{m_1}}=2^{\frac{1}{m}}$$
   and the reduction with respect to the Gauss valuation for $T=\rho_0X$ gives equation (1).
    \item Let's apply Proposition 3 a second time to the case where $r=\frac{m_1+1}{4}$, $s=\frac{m_2-m_1}{2}$, $(\rho_1,\rho_2)=(\varrho_1,\varrho_3)$  and $A_{\ell}$ lifting $b_{\ell}$. To do this, we easily check that the conditions $\frac{1}{2}v(2)\leq v(\varrho_1)\leq v(\varrho_3)$ et $v(\varrho_1^2\varrho_3)<2v(2)$ are satisfied. This gives  $n=2r+s=\frac{m_1+1}{2}+\frac{m_2-m_1}{2}=\frac{m_2+1}{2}$ and $\rho=\varrho_1^2\varrho_3$. We then find a cover of equation $Y^2=f_2(X)$ where $f_2(X)=\prod_{i=1}^{\frac{m_1+1}{4}}\tilde{P}_{i,2}(X)\prod_{j=r+1}^{r+s}P_{j,2}(X)$ and
    $$\tilde{P}_{i,2}(X)=\left[ (X-X_i)^2+\varrho_1\delta'_{i,2}(X-X_i)\right]^2+\varrho_1^2\varrho_3\beta'_{i,2}(X-X_i)^3+\varrho_1^2\varrho_3\alpha_{i,2}(X-X_i).
    $$
    $$P_{j,2}(X)=(X-X_j)^2+\varrho_1^2\varrho_3\gamma_{j,2}(X-X_j).$$
    In this case we have 
    $$ \rho_0=\left(\frac{4}{\rho}\right)^{\frac{1}{4r+2s-1}}=
    \left(\frac{4}{\varrho_1^2\varrho_3}\right)^{\frac{1}{m_2}}=2^{\frac{1}{m}}.$$
    and the reduction with respect to the Gauss valuation for $T=\rho_0X$ gives equation (2).
    \item Finally, let's apply Proposition 3 to the case where $r=\frac{m_1+1}{4}$, $s=\frac{m_3-m_1}{2}$ , $(\rho_1,\rho_2)=(\varrho_1,\varrho_2)$  and $A_{\ell}$ lifting $c_{\ell}$. To do this, we easily check that the conditions $\frac{1}{2}v(2)\leq v(\varrho_1)\leq v(\varrho_2)$ and $v(\varrho_1^2\varrho_2)<2v(2)$ are satisfied. Thus  $n=2r+s=\frac{m_1+1}{2}+\frac{m_3-m_1}{2}=\frac{m_3+1}{2}$ and $\rho=\varrho_1^2\varrho_2$. We then find a cover of equation  $Y^2=f_3(X)$ where $f_3(X)=\prod_{i=1}^{\frac{m_1+1}{4}}\tilde{P}_{i,3}(X)\prod_{j=r+1}^{r+s}P_{j,3}(X)$ and
    $$\tilde{P}_{i,3}(X)=\left[ (X-X_i)^2+\varrho_1\delta'_{i,3}(X-X_i)\right]^2+\varrho_1^2\varrho_2\beta'_{i,3}(X-X_i)^3+\varrho_1^2\varrho_2\alpha_{i,3}(X-X_i).
    $$
    $$P_{j,3}(X)=(X-X_j)^2+\varrho_1^2\varrho_2\gamma_{j,3}(X-X_j).$$
    In this case we have 
    $$ \rho_0=\left(\frac{4}{\rho}\right)^{\frac{1}{4r+2s-1}}=
    \left(\frac{4}{\varrho_1^2\varrho_2}\right)^{\frac{1}{m_3}}=2^{\frac{1}{m}}.$$
    and the reduction with respect to the Gauss valuation for $T=\rho_0X$ gives equation (3).
\end{enumerate}

%\subsection{Why it doesn't work yet}
Let's take stock of what's been proved.

The three covers thus constructed have simultaneously good reduction. Moreover, $Rev_1$, $Rev_2$ and $Rev_3$ have the points $X_1,\cdots, X_{\frac{m_1+1}{4}}$ as common branch points (which is the expected number as stated in Theorem 2). But Theorem 2 imposes additional combinatorial conditions : the coverings, considered in pairs, must have the right number of common branch points. The idea is therefore to modify the polynomials $\tilde{P}_{i,1},\tilde{P}_{i,2},\tilde{P}_{i,3}$ into new poynomials $\hat{P}_{i,1},\hat{P}_{i,2},\hat{P}_{i,3}$ such that 
\begin{enumerate}[a.]
    \item$X_i$ is always a root of $\hat{P}_{i,1},\hat{P}_{i,2},\hat{P}_{i,3}$ (to maintain the right number of branch points between the three covers).
    \item $\hat{P}_{i,1}$ and $\hat{P}_{i,2}$ have two common roots (one of which is already $X_i$). The same applies to $\hat{P}_{i,1}$ and $\hat{P}_{i,3}$, as well as to $\hat{P}_{i,2}$ and $\hat{P}_{i,3}$.
    \item The polynomials $\hat{P}_{i,1}$ and $\tilde{P}_{i,1}$ are "close enough", i.e. $\hat{P}_{i,1}=\tilde{P}_{i,1}+O(4)$ where $O(4)\in4R[X]$ is of degree at most 3. The same applies to $\hat{P}_{i,2}$ and $\tilde{P}_{i,2}$ as well as to $\hat{P}_{i,3}$ and $\tilde{P}_{i,3}$.
\end{enumerate}

If these conditions are met, the new polynomials (denoted $\hat{f}_1$, $\hat{f}_2$, $\hat{f}_3$) will have the right number of roots in common and we'll preserve the good simultaneous reduction of the covers.

\subsection{The key: Proposition 1 in ~\cite{Mat}}

Let's start with an easy lemma
\begin{lemma}
    Let $a_1,a_2\in R$. Then
    $$X(X+a_1^2)(X+a_2^2)(X+(a_1+a_2)^2)=(X^2+(a_1^2+a_2^2+a_1a_2)X)^2+a_1^2a_2^2(a_1+a_2)^2X.$$
\end{lemma}

\begin{Proof}\rm
    It's a special case of Proposition 1 of ~\cite{Mat} . The result stated there is true for all $n$ and is an equality modulo 4. It turns out that in the case $n=2$, the equality is exact. The case $n=2$ is also already mentioned in ~\cite{Mat3}.

\begin{flushright}
     $\square$
     \end{flushright}
     
    Having this result in reserve for all $n$ bodes well for an eventual generalization to the $(\Z/2\Z)^n$ case.

\end{Proof}
Let's now consider, for any $i\in\CG 1,\frac{m_1+1}{4}\CD$, elements $a_{i,1}, a_{i,2},a_{i,3}\in R$. Let
$$
\begin{aligned}
    \PP_{i,1}(X)&=(X-X_i)(X-X_i+a_{i,2}^2)(X-X_i+a_{i,3}^2)(X-X_i+(a_{i,2}+a_{i,3})^2)\\
    \PP_{i,2}(X)&=(X-X_i)(X-X_i+a_{i,1}^2)(X-X_i+a_{i,3}^2)(X-X_i+(a_{i,1}+a_{i,3})^2)\\
    \PP_{i,3}(X)&=(X-X_i)(X-X_i+a_{i,1}^2)(X-X_i+a_{i,2}^2)(X-X_i+(a_{i,1}+a_{i,2})^2)
\end{aligned}$$

Note that these three polynomials have the "right combinatorics" in terms of their roots. To get the final result, we'd have to modify them so that they coincide (modulo 4) with the polynomials $\tilde{P}_{i,1},\tilde{P}_{i,2},\tilde{P}_{i,3}$.

In particular, Lemma 5 tells us that 
$$
\begin{aligned}
    \PP_{i,1}(X)&=((X-X_i)^2+(a_{i,2}^2+a_{i,3}^2+a_{i,2}a_{i,3})(X-X_i))^2+a_{i,2}^2a_{i,3}^2(a_{i,2}+a_{i,3})^2(X-X_i)\\
    \PP_{i,2}(X)&=((X-X_i)^2+(a_{i,1}^2+a_{i,3}^2+a_{i,1}a_{i,3})(X-X_i))^2+a_{i,1}^2a_{i,3}^2(a_{i,1}+a_{i,3})^2(X-X_i)\\
    \PP_{i,3}(X)&=((X-X_i)^2+(a_{i,1}^2+a_{i,2}^2+a_{i,1}a_{i,2})(X-X_i))^2+a_{i,1}^2a_{i,2}^2(a_{i,1}+a_{i,2})^2(X-X_i)
\end{aligned}$$

It would therefore be necessary to
\begin{enumerate}[a.]
    \item Choose $a_{i,1}, a_{i,2},a_{i,3}\in R$ such that
    $$\left\{
    \begin{aligned}
        a_{i,2}^2a_{i,3}^2(a_{i,2}+a_{i,3})^2&=\varrho_2^2\varrho_3\alpha_{i,1}+O(4)\\
        a_{i,1}^2a_{i,3}^2(a_{i,1}+a_{i,3})^2&=\varrho_1^2\varrho_3 \alpha_{i,2}+O(4)\\
        a_{i,1}^2a_{i,2}^2(a_{i,1}+a_{i,2})^2&=\varrho_1^2\varrho_2 \alpha_{i,3}+O(4)
    \end{aligned}
    \right.$$
    \item Make appear terms of the form $\varrho_2^2\varrho_3\beta'_{i,1}(X-X_i)^3$ (resp. $\varrho_1^2\varrho_3\beta'_{i,2}(X-X_i)^3$, $\varrho_1^2\varrho_2\beta'_{i,3}(X-X_i)^3$) in $\PP_{i,1}(X)$ (resp. $\PP_{i,2}(X)$, $\PP_{i,3}(X)$).
    \item Choose the $\delta'_{i,1}, \delta'_{i,2}, \delta'_{i,3}$ as required, but this will be easy. 
\end{enumerate}
\subsection{End of the proof of Theorem 1.b}

The most complicated point is point a. This is the subject of the following proposition :
\begin{Prop}
   With the previous notations, we can find $a_{i,1}, a_{i,2},a_{i,3}\in R$ such that 
    $$\left\{
    \begin{aligned}
        a_{i,2}^2a_{i,3}^2(a_{i,2}+a_{i,3})^2&=\varrho_2^2\varrho_3\alpha_{i,1}+O(4)\\
        a_{i,1}^2a_{i,3}^2(a_{i,1}+a_{i,3})^2&=\varrho_1^2\varrho_3 \alpha_{i,2}+O(4)\;\;\;\;(S)\\
        a_{i,1}^2a_{i,2}^2(a_{i,1}+a_{i,2})^2&=\varrho_1^2\varrho_2 \alpha_{i,3}+O(4)
    \end{aligned}
    \right.$$
    More precisely, if we write $u_1=(\varrho_2^2\varrho_3\alpha_{i,1})^{\frac{1}{4}}$,
     $u_2=(\varrho_1^2\varrho_3\alpha_{i,2})^{\frac{1}{4}}$ et $u_3=(\varrho_1^2\varrho_2\alpha_{i,3})^{\frac{1}{4}}$, just take
     $$\left\{
        \begin{aligned}
            a_{i,1}&=\frac{u_2u_3(u_2+u_3)}{(u_1u_2u_3(u_1+u_2)(u_1+u_3)(u_2+u_3)(u_1+u_2+u_3))^{\frac{1}{3}}}\\
            a_{i,2}&=\frac{u_1u_3(u_1+u_3)}{(u_1u_2u_3(u_1+u_2)(u_1+u_3)(u_2+u_3)(u_1+u_2+u_3))^{\frac{1}{3}}}\\
           a_{i,3}&=\frac{u_1u_2(u_1+u_2)}{(u_1u_2u_3(u_1+u_2)(u_1+u_3)(u_2+u_3)(u_1+u_2+u_3))^{\frac{1}{3}}} 
        \end{aligned}
     \right.$$

     Furthermore, $v(a_{i,j})=\frac{1}{2}v(\varrho_j)$ for all $i\in\cg 1,r\cd]$ and $j\in\cg1,3\cd$.
\end{Prop}

\begin{Remark}\rm
   We can use the notations from ~\cite{Mat} (page 96). We would then write
    
 $$a_{i,1}=\frac{\pi_2(u_2,u_3)}{\pi_3(u_1,u_2,u_3)^{\frac{1}{3}}},\;\;
 a_{i,2}=\frac{\pi_2(u_1,u_3)}{\pi_3(u_1,u_2,u_3)^{\frac{1}{3}}},\;\;
 a_{i,3}=\frac{\pi_2(u_1,u_2)}{\pi_3(u_1,u_2,u_3)^{\frac{1}{3}}}\;\;$$
\end{Remark}

\begin{Proof}\rm
     The proof is divided into three parts :
      \begin{enumerate}[i.]
          \item Justify that $v(u_i+u_j)=\min(v(u_i),v(u_j))$ as soon as $i\neq j$ (and also that $v(u_1+u_2+u_3)=\min(v(u_1),v(u_2),v(u_3))$.
          \item Verify that the $a_{i,j}$ defined by our formulas satisfy the relation $v(a_{i,j}^2)=v(\varrho_j)$ for $j\in\cg1,3\cd$.
          \item Check that the $a_{i,j}$ thus defined satisfy the equations of the system $(S)$. 
      \end{enumerate}

      \begin{enumerate}[i.]
          \item First, let's note that the elements $\alpha_{i,1}, \alpha_{i,2}, \alpha_{i,3}$ have valuation zero (we refer to the definition of $\alpha_i$ at the beginning of part V.3). 
          
          In the case where $m_1<m_2<m_3$, we have $v(\varrho_1)<v(\varrho_2)<v(\varrho_3)$ so that $u_1$, $u_2$ and $u_3$ have different valuations. The result is immediate. 

          In the case where $m_1=m_2<m_3$, we have $v(\varrho_1)=v(\varrho_2)<v(\varrho_3)$, so that $u_1$ and $u_2$ have the same valuation. We must then justify that $u_1+u_2$ has the same valuation as $u_1$ and $u_2$. But
          $$u_1+u_2=(\varrho_1^2\varrho_3)^{\frac{1}{4}}(\alpha_{i,1}^{\frac{1}{4}}+\alpha_{i,2}^{\frac{1}{4}})$$
          so we need to check that $(\alpha_{i,1}^{\frac{1}{4}}+\alpha_{i,2}^{\frac{1}{4}})$ is not zero in reduction modulo $\pi$, i.e. that $(\alpha_{i,1}+\alpha_{i,2})$ is not zero in reduction modulo $\pi$.

          We have $\alpha_{i,1}=\sum_{\ell=0}^{n-1}\rho_0^{2\ell}A_{\ell}^2\alpha_{i\ell}^2$ where $A_{\ell}$ lift $a_{\ell}$ and $\alpha_{i\ell}^2$ are the coefficients of the partial fraction expansion of $\frac{X^{2\ell}}{\prod\limits_{i=1}^r(X-X_i)^4}$. The same goes for $\alpha_{i,2}=\sum_{\ell=0}^{n-1}\rho_0^{2\ell}A_{\ell}^2\alpha_{i\ell}^2$ where this time $A_{\ell}$ lift $b_{\ell}$ and the $\alpha_{i\ell}^2$ are always  the coefficients of the partial fraction expansion of $\frac{X^{2\ell}}{\prod\limits_{i=1}^r(X-X_i)^4}$ (it's the same fraction  since $m_1=m_2$ and then $s=0$). The reduction modulo $\pi$ of $\alpha_{i,1}+\alpha_{i,2}$ is therefore $(a_0^2+b_0^2)\overline{\alpha_{0\ell}}$ and this element is non-zero.

         The cases $m_1<m_2=m_3$ and $m_1=m_2=m_3$ are treated in the same way. 
        \item Considering what has just been shown, we know the valuations of all the terms involved in the expression of $a_{i,j}$. Let's do the calculation for $a_{i,1}$ for example. We have
        $$\begin{aligned}
            v(a_{i,1})&=2v(u_3)+v(u_2)-\frac{1}{3}(4v(u_3)+2v(u_2)+v(u_1))\\
            &=\frac{1}{3}(2v(u_3)+v(u_2)-v(u_1))\\
            &=\frac{1}{3}(v(\varrho_1)+\frac{1}{2}v(\varrho_2)+\frac{1}{2}v(\varrho_1)+\frac{1}{4}v(\varrho_3)-\frac{1}{2}v(\varrho_2)-\frac{1}{4}v(\varrho_3))\\
            &=\frac{1}{2}v(\varrho_1)
        \end{aligned}$$
       The calculation is similar for $a_{i,2}$ and $a_{i,3}$. 
       \item Let's check the last equation in the system (the other two are similar). We have :
       $$\begin{aligned}
           a_{i,1}a_{i,2}(a_{i,1}+a_{i,2})&=\frac{u_1u_2u_3^3(u_1+u_3)(u_2+u_3)(u_1^2+u_2^2+u_3(u_1+u_2))}{u_1u_2u_3(u_1+u_2)(u_1+u_3)(u_2+u_3)(u_1+u_2+u_3)}\\
           &=\frac{u_3^2(u_1^2+u_2^2+u_3(u_1+u_2))}{(u_1+u_2)(u_1+u_2+u_3)}\\
           &=u_3^2\frac{(u_1+u_2)(u_1+u_2+u_3)-2u_1u_2}{(u_1+u_2)(u_1+u_2+u_3)}\\
           &=u_3^2-\frac{2u_1u_2u_3^2}{(u_1+u_2)(u_1+u_2+u_3)}\\
           &=u_3^2-2N
       \end{aligned}$$
by writing $N=\frac{u_1u_2u_3^2}{(u_1+u_2)(u_1+u_2+u_3)}$. It's easy to see that $v(N)\geq 0$ and therefore that $N\in R$.

    Therefore
       $$ a_{i,1}^2a_{i,2}^2(a_{i,1}+a_{i,2})^2=u_3^4+4N^2-4Nu_3^2=\varrho_1^2\varrho_2 \alpha_{i,3}+O(4).$$
      \end{enumerate}

      \begin{flushright}
     $\square$
     \end{flushright}
\end{Proof}

Let's make this choice for the elements $a_{i,1}$, $a_{i,2}$ and $a_{i,3}$ and modify the polynomials $\PP_{i,1},\PP_{i,1},\PP_{i,1}$ into new polynomials $\hat{P}_{i,1}, \hat{P}_{i,2},\hat{P}_{i,3}$ defined by 
$$
\begin{aligned}
    \hat{P}_{i,1}(X)&=(X-X_i)(X-X_i+a_{i,2}^2)(X-X_i+a_{i,3}^2)(X-X_i+(a_{i,2}+a_{i,3})^2+\varrho_2^2\varrho_3\beta'_{i,1}))\\
    \hat{P}_{i,2}(X)&=(X-X_i)(X-X_i+a_{i,1}^2)(X-X_i+a_{i,3}^2)(X-X_i+(a_{i,1}+a_{i,3})^2+\varrho_1^2\varrho_3\beta'_{i,2}))\\
    \hat{P}_{i,3}(X)&=(X-X_i)(X-X_i+a_{i,1}^2)(X-X_i+a_{i,2}^2)(X-X_i+(a_{i,1}+a_{i,2})^2+\varrho_1^2\varrho_2\beta'_{i,3}))
\end{aligned}$$
Recall that the elements $\beta_i'$ were defined in section V.4 to state Proposition 3.
Let's look at the first of these polynomials. We have
$$\begin{aligned}
    \hat{P}_{i,1}(X)&=\PP_{i,1}(X)+\varrho_2^2\varrho_3\beta'_{i,1}(X-X_i)(X-X_i+a_{i,2}^2)(X-X_i+a_{i,3}^2)\\&=\PP_{i,1}(X)+\varrho_2^2\varrho_3\beta'_{i,1}(X-X_i)^3+O(4)\\
    &=((X-X_i)^2+(a_{i,2}^2+a_{i,3}^2+a_{i,2}a_{i,3})(X-X_i))^2\\
    &+\varrho_2^2\varrho_3\alpha_{i,1}(X-X_i)+\varrho_2^2\varrho_3\beta'_{i,1}(X-X_i)^3+O(4)
\end{aligned}$$

It remains to choose $\varrho_2\delta'_{i,1}$ equal to $(a_{i,2}^2+a_{i,3}^2+a_{i,2}a_{i,3})$, to be sure that $\hat{P}_{i,1}(X)$ verifies the hypothesis of Proposition 3. This is possible since we can be sure that the valuations of each of the terms $a_{i,2}^2,a_{i,3}^2$ and $a_{i,2}a_{i,3}$ are greater than or equal to that of $\varrho_2$.
The operation is the same with $ \hat{P}_{i,2}(X)$ and $ \hat{P}_{i,3}(X)$. 

It remains to check that the polynomials $\hat{P}_{i,1},\hat{P}_{i,2},\hat{P}_{i,3}$, don't have too many roots in common, i.e. that the seven elements
 $$0,a_{i,1}^2,a_{i,2}^2,a_{i,3}^2,(a_{i,2}+a_{i,3})^2+\varrho_2^2\varrho_3\beta'_{i,1},(a_{i,1}+a_{i,3})^2+\varrho_1^2\varrho_3\beta'_{i,2},(a_{i,1}+a_{i,2})^2+\varrho_1^2\varrho_2\beta'_{i,3}$$
are indeed distinct in pairs. If this is not the case, we can modify one or more of these seven elements by adding an element of $4R$, to make them distinct.

\begin{Remark} Probably these seven elements are automatically distinct, but it's a bit long to check, as one would have to separate cases, depending on whether the conductors are equal or not. The previous argument avoids the need for painful justifications.
    
\end{Remark}

So the three covers given by the equations $Y^2= \hat{P}_{i,1}(X)$ have simultaneously good reduction and the right number of branch points in common, which completes
the proof of Theorem 1.b, and thus of the theorem as a whole.
\section{The case $(\Z/2\Z)^n$}

  We can hope to generalize this result to $(\Z/2\Z)^n$ with $n\geq 4$, with conductors $m_1+1\leq\cdots\leq m_n+1$ (this $n$-tuple always being assumed minimal) verifying the only divisibility condition: $2^{n-i}|m_i+1$ for $i\in\cg1,n-1\cd$ and $2|m_n+1$.

   Finding the " good " geometry of the branch locus with the right thicknesses is not a problem. It also seems that the construction of the covers with imposed geometry and reduction  (as in Part IV) can be done in the same way. The remaining difficulty is the combinatorics of the branch locus.

%\documentclass{article}
%\usepackage{graphicx} % Required for inserting images

%\title{(Zsur2Z)3-lifting}
%\author{Guillaume Pagot}
%\date{August 2023}

%\begin{document}

%\maketitle

%\section{Introduction}

\end{document}